\documentclass[graybox]{svmult}

\usepackage{type1cm}        % activate if the above 3 fonts are
\usepackage{makeidx}         % allows index generation
\usepackage{graphicx}        % standard LaTeX graphics tool
\usepackage{multicol}        % used for the two-column index
\usepackage[bottom]{footmisc}% places footnotes at page bottom

\usepackage{newtxtext}       % 
\usepackage{newtxmath}       % selects Times Roman as basic font

\usepackage{amsmath,amsfonts,mathrsfs,paralist}
\usepackage{mathtools,enumitem}
\usepackage{setspace}
\usepackage{dsfont}
\usepackage{tcolorbox}

\newcommand{\N}{\mathbb {N}}

\def\BFu {\mathbf{u}}
\def\BFe {\mathbf{e}}

\def\BFU {\mathbf{U}}
\def\BFV {\mathbf{V}}
\def\BFv {\mathbf{v}}
\def\BFx {\mathbf{x}}

\def\BFy {\mathbf{y}}

\def\BFi {\textit{\bfseries i}}

\def\BFk {\textit{\bfseries k}}
\def\BFd {\textit{\bfseries d}}

\def\Ik {\mathbf{1}_{\BFk}}

\def\iddots{\mathinner{\mkern1mu\raise1pt
    \hbox{.}\mkern2mu\raise4pt\hbox{.}\mkern2mu
        \raise7pt\vbox{\kern7pt\hbox{.}}\mkern1mu}}
\def\Vol{{\rm Vol}}

\DeclareMathOperator{\vol}{Vol}

\def\hjaspar#1{}
\def\hchristiane#1 {}

\makeindex             % used for the subject index
\begin{document}

\title*{On the distribution of scrambled $(0,m,s)$-nets over unanchored boxes}
\author{C. Lemieux and J. Wiart}
\institute{C. Lemieux \at University of Waterloo, 200 University Ave.~West ON, Canada, N2L 3G1, \email{clemieux@uwaterloo.ca}
\and J.Wiart \at Johannes Kepler University, Altenbergerstr.\ 69,
4040 Linz, Austria, \email{jaspar.wiart@jku.at}}
\maketitle

%Each chapter should be preceded by an abstract (no more than 200 words) that summarizes the content. The abstract will appear \textit{online} at \url{www.SpringerLink.com} and be available with unrestricted access. This allows unregistered users to read the abstract as a teaser for the complete chapter.
%Please use the 'starred' version of the \texttt{abstract} command for typesetting the text of the online abstracts (cf. source file of this chapter template \texttt{abstract}) and include them with the source files of your manuscript. Use the plain \texttt{abstract} command if the abstract is also to appear in the printed version of the book.}

\abstract{We introduce a new quality measure to assess randomized low-discrepancy point sets of finite size $n$. This new quality measure, which we call ``pairwise sampling dependence index'', is based on the concept of negative dependence. A negative value for this index implies that the corresponding point set integrates the indicator function of any unanchored box with smaller variance than the Monte Carlo method.
We show that scrambled $(0,m,s)-$nets have a negative pairwise sampling dependence index. We also illustrate through an example that randomizing via a digital shift instead of scrambling may yield a positive pairwise sampling dependence index.}
%Each chapter should be preceded by an abstract (no more than 200 words) that summarizes the content. The abstract will appear \textit{online} at \url{www.SpringerLink.com} and be available with unrestricted access. This allows unregistered users to read the abstract as a teaser for the complete chapter.\newline\indent
%Please use the 'starred' version of the \texttt{abstract} command for typesetting the text of the online abstracts (cf. source file of this chapter template \texttt{abstract}) and include them with the source files of your manuscript. Use the plain \texttt{abstract} command if the abstract is also to appear in the printed version of the book.}

\hchristiane{Check $u_{ij}$ vs $u_{i,j}$ (only saw it in Def 1) and $l$ vs $\ell$ (I think I fixed those as well)}
\section{Introduction}
\label{sec:intro}

The quality of point sets used within quasi-Monte Carlo (QMC) methods is often assessed using the notion of discrepancy. For a point set $P_n=\{\BFu_i:i=1,\ldots,n\}$, its star-discrepancy is given by
$
D^*_n(P_n) = \sup_{A \in {\cal A}_0} |J_n(A)-{\rm Vol(A)}| 
$
where ${\cal A}_0$ is the set of all boxes $A \subseteq [0,1)^s$ anchored at the origin, and $J_n(A) = \sum_{i=1}^n \mathbf{1}_{\BFu_i \in A}/n$. The extreme discrepancy is instead given by 
$
D_n(P_n) = \sup_{A \in {\cal A}} |J_n(A)-{\rm Vol(A)}| 
$
where ${\cal A}$ is the set of all boxes in $[0,1)^s$. Both quantities are typically interpreted as comparing the empirical distribution induced by $P_n$ with the uniform distribution over $[0,1)^s$ in terms of the probability they assign to a given set ${\cal A}$ of boxes. Using inclusion-exclusion arguments, one can derive the bound $D_n(P_n) \le 2^s D^*_n(P_n)$.

Many asymptotic results for $D^*_n(P_n)$ and $D_n(P_n)$ have been derived for various low-discrepancy sequences \cite{DiPi10,rNIE92b}. These sequences are understood to be such that $D^*_n(P_n) \in O((\log n)^s/n)$, and the above mentioned results often focus on studying the constant terms in the big-Oh notation and how it behaves as a function of $s$. 

In practice, when using QMC methods, one is often working in settings where $n$ is not too large, and a primary goal is to make sure that the QMC approximation will result in a better approximation than the one that would be obtained by using plain Monte Carlo sampling. One is also typically interested in assessing the approximation error, something naturally embedded in Monte Carlo methods via variance estimates and the central limit theorem.

In this setting, the use of randomized QMC methods is very appealing, as it preserves the advantage of QMC induced by the use of low-discrepancy sequences, while at the same time allowing for error estimates through independent and identically distributed (iid) replications. 

In this paper, we focus on the above settings, i.e., where one  (1) works with $n$ not too large; (2) uses randomized QMC, and (3) hopes to do better than Monte Carlo.

\hchristiane{state main result and give sketch of proof}

To this end, we propose to reinterpret the measures $D^*_n(P_n)$ and $D_n(P_n)$ and propose a new, related measure that is designed for our chosen setting, which we refer to as ``pairwise sampling dependence index''. While this measure is meant to assess the uniformity of point sets much like the star and extreme discrepancies do, it also has another interpretation, which is that a point set with negative pairwise sampling dependence index estimates the expected value of the indicator function $\mathbf{1}_A$ for any $A \in {\cal A}$ with variance smaller than the Monte Carlo method. This new measure is  defined in Section \ref{sec:disc},  Eqn.\ \eqref{eq:pairwise sdi}. In Section \ref{sec:anchored} we revisit the result from \cite{WLD20}, which shows that scrambled  \hjaspar{I believe this is only for digital nets} $(0,m,s)-$nets have a negative pairwise sampling dependence index over all anchored boxes. %if and only if $t=0$. 
In Section \ref{sec:negdepunanchoredbox} we show that this result extends to unanchored boxes {in  Theorem \ref{thm:scr0msnegdepunanchored}, which is the main result of this paper.} Hence the extension to unanchored boxes does not cause the same deterioration of the bound for this uniformity measure  as is the case when applying an inclusion-exclusion argument to go from the star to the extreme discrepancy. %The proof of this result makes use of a carefully constructed decomposition of unanchored boxes that is presented in Section \ref{sec:decomp}.

{The proof of this result is essentially a very difficult problem in linear programming, (something that is rather obscured by the fact that, since the number of variables depends on $m$, we work in $\ell^1(\N)$ and its dual $\ell^\infty(\N)$ rather than a finite dimensional space). Indeed, we must demonstrate that \eqref{eq:pairwise sdi} is always negative for scrambled $(0,m,s)$-nets. We see in  Theorem \ref{thm:decompHA} \hjaspar{By moving the definition of $H(A)$ to \eqref{eq:pairwise sdi}, this requires less hunting for definitions.} that \eqref{eq:pairwise sdi} is actually a linear equation whose variables are non-negative and are further constrained, in the one-dimensional case, according to Lemma \ref{lem:vol}. The constraints define a convex region whose extreme points, in the one-dimensional case, are found in Theorem \ref{thm:decompA} and given by \eqref{eq:ExtreemPoints} in the higher dimensional case. The remainder of the proof boils down to proving that \eqref{eq:pairwise sdi} is non-negative at these extreme points.} This requires the use of
%It also uses 
several technical lemmas (given in the appendix) proving sufficiently tight bounds on various combinatorial sums, which is precisely why we do not have to rely on an inclusion-exclusion argument to go from the anchored case to the unanchored one. We briefly discuss in Section \ref{sec:scrambling} the advantage of scrambling over simpler randomization methods such as a digital shift. Ideas for future work are presented in Section \ref{sec:conc}.

%A version of this paper that includes the technical lemmas and proofs required to get our main results %from Section \ref{sec:negdepunanchoredbox} 
%can be found on the arXiv.

\section{Pairwise Sampling Dependence}
\label{sec:disc}

We start by revisiting the definition of extreme discrepancy using a probabilistic approach, despite the fact that the point set $P_n$ may be deterministic. We do so by introducing a quality measure we call {\em sampling discrepancy}, given by 
$
{\cal D}_n (P_n) := \sup_{A \in {\cal A}}|{\cal P}_n(A) - {\rm Vol}(A)|,
$
where ${\cal P}_n(A)$ is the probability that a randomly chosen point in $P_n$ will fall in $A$.
%, and $\lambda(\cdot)$ is the Lebesgue measure. 
For a deterministic point set, this probability is given by ${\cal P}_n(A) = J_n(A)/n$, and so in this case
${\cal D}_n(P_n) = D_n(P_n)$. This definition captures how the discrepancy is often described as a distance measure between the empirical distribution induced by the point set $P_n$ and the uniform distribution. Since the uniform distribution is viewed as a \emph{target distribution} in this setting, we want this distance to be as small as possible.

In this paper we are interested in randomized QMC point sets $\tilde{P}_n$. We assume  $\tilde{P}_n$ is a valid sampling scheme, meaning that $\BFU_i \sim U(0,1)^s$ for each $\BFU_i \in \tilde{P}_n$, with possibly some dependence among the $\BFU_i$'s. When we write $\tilde{P}_n$, we are thus not referring to a specific realization of the randomization process, which we instead denote by $\tilde{P}_n(\omega)$, where $\omega \in \Omega$, and $\Omega$ is the sampling space associated with our randomization process for $P_n$.

In that setting, we could compute ${\cal D}_n(\tilde{P}_n(\omega))$ %for a given realization $P_n(\omega)$, 
and then perhaps compute the expected value of  ${\cal D}_n(\tilde{P}_n(\omega))$ over all these realizations $\omega$, or the probability that it will be larger than some value, as done in \cite{Gnewuch2019}, for example. 
If we instead interpret ${\cal P}_n(A)$ as the probability that a randomly chosen point $\BFU_i$ from $\tilde{P}_n$ falls in $A$, then we would get ${\cal D}_n(\tilde{P}_n)=0$, which is of little use.

%Typically, when studying the discrepancy of randomized point sets, one avoids this problem by inserting the expectation operator outside the
%distance measure,  e.g., as in the definition of the mean square $L^2$-discrepancy, given by 
%\[
%{\rm E} \int_{\BFx} (J_n(\BFx) - {\rm %Vol}([0,\BFx))^2 d\BFx
%
To define an interesting alternative measure of uniformity for randomized QMC point sets,  we introduce instead a ``second-moment'' version of the sampling discrepancy, in which we consider pairs of points rather than single points, and where the distribution against which we compare the point set is that induced by random sampling, where points are sampled independently from one another. When considering pairs of points, our goal is to examine the propensity for points to repel each other, which is a desirable feature if we want to achieve greater uniformity than random sampling. Note that this notion of ``propensity to repel'' is in line with the concept of negative dependence. 

More precisely, we want to verify that the pairs of points from $\tilde{P}_n$ are less likely to fall within the same box $A$ than they would if they were independent. %, as is the case with random sampling.
Note that here, as was the case with ${\cal D}_n(P_n)$, we are comparing the distribution induced by $\tilde{P}_n$ with another distribution. But rather than comparing to a target distribution to which we want to be as close as possible, we are comparing to a distribution upon which we want to improve, and thus are not trying to be close to that distribution.

%So this is similarly to what is done with pointwise sampling discrepancy, we are comparing the distribution induced by the point set with another distribution, but here, we compare to random sampling rather than the uniform distribution (we need a distribution that tells us how pairs of points behave). Furthermore, we are not trying to mimic that target distribution, we are instead trying to show that rather than trying to be as close as possible to t 
The measure we propose  to assess the quality of a point set via the behavior of its pairs  is called \emph{pairwise sampling dependence index}  and is given by \hjaspar{Consider moving the definition of $(A)$ here}
\begin{align}
{\cal E}_n(\tilde{P}_n) &:= \sup_{A \in {\cal A}} H_n(A)  - {\rm Vol}^2(A),\label{eq:pairwise sdi}\\
\mbox{where }    H_n(A) &:= P((\BFU,\BFV) \in A \times A),
\label{eqdef:HnA}
\end{align}
with $\BFU$ and $\BFV$ being distinct points randomly chosen from $\tilde{P}_n$. %We call this quantity 
%${\cal E}_n(\tilde{P}_n)$ 
We say $\tilde{P}_n$ has a \emph{negative pairwise sampling dependence index} when ${\cal E}_n(\tilde{P}_n) \le 0.$
(This terminology is consistent with other measures of negative dependence: see, e.g.,  \cite{WLD20}.)

%In other words, here we consider the probability that a pair of points randomly chosen from $P_n$ will be in $A \times A$, and evaluate the difference between the distribution induced by $P_n$ and random sampling. 
Note that because we are not taking the supremum over all products of the form $A \times B$ with $A,B \in {\cal A}$ and instead only consider $A \times A$, it is possible, if $\tilde{P}_n$ is designed so that points tend to cluster away from each other, that the probability $H_n(A)$ %$P((\BFU,\BFV) \in A \times A)$ 
will never be larger than what it is under random sampling, as given by ${\rm Vol}^2(A)$.
This is not the case with the measure ${\cal D}_n(P_n)$, where having ${\cal P}_n(A)< {\rm Vol}(A)$ implies there will be some $A'$ for which ${\cal P}_n(A')> 
{\rm Vol}(A')$.
There the goal is to show there exist point sets $P_n$ with $|{\cal P}_n(A)-{\rm Vol}(A)|$ very close to 0, and becoming closer to 0 as $n$ goes to infinity. In our case, we instead want to show, for a given $n$, that there exist sampling schemes $\tilde{P}_n$ with ${\cal E}_n(\tilde{P}_n) \le 0$.

{So far we mentioned the work done in \cite{Gnewuch2019} and \cite{WLD20}, but concepts of dependence based on measures different from \eqref{eq:pairwise sdi} have recently been used in other works to analyze lattices \cite{Wnuk19lat}, Latin hypercube sampling \cite{Gnewuch2021a} and scrambled nets \cite{Doerr2021}.  }
%For example, \christiane{Need to add references somewhere to Wnuk and Gnewuch on lattices (ORL), Gnewuch and Hebbinghaus (... LHS, arXiv) and Doerr and Gnewuch (neg dep properties of LHS and scrambled nets) and explain connection with our work}}

%%%%%%%%%%%%%%%%%%%%%%%%%%%%%%%%%%%%%
\section{Revisiting Pairwise Sampling Dependence over Anchored Boxes}
\label{sec:anchored}

In what follows, we assume $\tilde{P}_n$ is a scrambled $(0,m,s)-$net in base $b \ge s$, where $n=b^m$, and $P_n$ represents the underlying $(0,m,s)-$net being scrambled.
We assume the reader is familiar with the concept of digital nets and $(t,m,s)-$nets, as presented in e.g., \cite{DiPi10,rNIE92b}. Also, when referring to scrambled nets, we refer to the scrambling method studied in \cite{WLD20}, which originates from \cite{vOWE95a}.

In \cite{WLD20}, it was shown that  if we restrict ${\cal E}_n(\tilde{P}_n)$ to anchored boxes% $A \in {\cal A}_0$
---denote this version of ${\cal E}_n$ by ${\cal E}_{n,0}$---then ${\cal E}_{n,0}(\tilde{P}_n) \le 0$. In fact, a stronger result is shown in  \cite{WLD20}, which is that  for $(\BFU,\BFV)$ a pair of distinct points randomly chosen from $\tilde{P}_n$, %the inequality
\[
P((\BFU,\BFV) \in [\mathbf{0},\BFx) \times [\mathbf{0},\BFy)) \le {\rm Vol}([\mathbf{0},\BFx) \times [\mathbf{0},\BFy)), \qquad \mbox{for any }
\BFx,\BFy \in [0,1]^s.
\]
For simplicity, in what follows we assume $\BFx = \BFy$ and let $A = [\mathbf{0},\BFx).$ 
%Furthermore, we let $H(A)$
%be defined as 
%\[
%H(A) := P((\BFU,\BFV) \in A \times A).
%\]
\hchristiane{Ref 2 \# 18 says maybe we need to use calligraphic $V$ instead of plain $V$}

{Next, to define a key quantity called the {\em volume vector} of a subset of $[0,1]^{2s}$, we first define the regions} $D_{\BFi}^s := \{(\BFx,\BFy) \in [0,1)^{2s}: \boldsymbol{\gamma}_b^s(\BFx,\BFy) = \BFi\}$, where
$	\boldsymbol{\gamma}^s_b(\BFx,\BFy) := (\gamma_b(x_1,y_1),\ldots,$ $\gamma_b(x_s,y_s))$
and $\gamma_b(x,y) \ge 0$ is 
the {unique} number $i \ge 0$ such that \hchristiane{ref 2 says to say that $i$ is unique: I don't think I agree}
	\begin{equation} \label{eq:gammadef}
	\lfloor b^ix\rfloor=\lfloor b^iy\rfloor\quad\text{but}\quad\lfloor b^{i+1}x\rfloor \neq\lfloor b^{i+1}y\rfloor.
	\end{equation}
That is, $\gamma_b(x,y)$ is the exact number of initial digits shared by $x$ and $y$ in their base $b$ expansion. 	If $x=y$ then we let $\gamma_b(x,y) = \infty$. Also, \eqref{eq:gammadef} implies  $\gamma_b(x,y)$ is well defined for any $x,y \in [0,1)$ even if  $x,y$ do not have a unique expansion in base $b$.

Let $\N_0 = \{0,1,2,\ldots\}$. We can now define, for $A,B \subseteq [0,1]^s$, the volume vector $V(A \times B) \in \ell^1(\N_0^s)$, whose component $V_{\BFi}(A \times B)$ associated to $\BFi \in \N_0^s$ \hchristiane{Ref 2 says that including 0 in $\N$ will mislead people but doesn't suggest what to do instead! Maybe we can rebuke by giving examples of other places where it is defined like this.} 
is given by
\hchristiane{Ref 2 doesn't like that we go in reverse order: would prefer to build up notation and motivate it}\hjaspar{The notation in the definition does not match the notation in equation \ref{eq:HAinsidelemma}}
\[
V_{\BFi}(A \times B) := \int_{A \times B} {\bf 1}_{D_{\BFi}^s} d\BFu d\BFv = \vol((A \times B)\cap D_{\BFi}^s) \in [0,1].
\]

A key step  used in \cite{WLD20} to prove that $H_n(A) \le {\rm Vol}(A \times A)$ is to find a conical combination of products of the form $\Ik \times \Ik$, where $\Ik = \prod_{j=1}^s [0,b^{-k_j})$, $k_j \in \N_0, j=1,\ldots,s,$ whose volume vector is the same as that of $A \times A$. More precisely, one can find coefficients $t_{\BFk} \ge 0$ with
$\sum_{\BFk \ge \mathbf{0}} t_{\BFk} = {\rm Vol}(A \times A)$,
such that
\hchristiane{Ref 2 wants us to add subscript $\BFi$ to $V(A)$ and $V(\Ik)$ and then say this holds for all $\BFi$: reasoning is that $V$ is much more complicated than it looks in that equation: check that it's equivalent}
\begin{equation}
\label{eq:VAdecompElemk}
{
V_{\BFi}(A \times A) = \sum_{\BFk\ge \mathbf{0}} t_{\BFk} b^{2|\BFk|} V_{\BFi}(\Ik \times \Ik) \qquad \mbox{for all } \BFi \in \N_0^s.}
\end{equation}
The coefficients $t_{\BFk}$ are shown in \cite{WLD20} to be given by $t_{\BFk}=\prod_{j=1}^s t_{k_j}$, where
\begin{equation}
    \label{eq:deftk}
t_{k} = \begin{cases}
\frac{bV_k(A \times A)-V_{k-1}(A \times A)}{b-1} & \mbox{ if } k>0 \\
\frac{bV_0(A \times A)}{b-1} & \mbox{ if } k=0.
\end{cases}
\end{equation}

%T See [Negative dependence paper, Section 2.3] for a more detailed explanation.

From here, rather than following the proof in \cite{WLD20}, we
{exploit the fact that the joint pdf of points $(\BFU,\BFV)$ from a scrambled $(0,m,s)-$net is a simple function that is constant on the $D_{\BFi}$ regions. The volume vector simply keeps track of how much the $D_{\BFi}$ region is covered by $A\times A$, allowing $H_n(A)$ to be written as a linear sum. Since the joint pdf is a simple function, these sums always have finitely many non-zero terms. See \cite[Sec.\ 2.3]{WLD20} for more details.}
%instead make use of
%a property of the volume vector operator, which gives the following %representation:
\hchristiane{Check new wording suggested by ref 2, i.e., let $ t_{\BFk}$ be the coefficients from (2) instead of saying ``for which (2) holds''}
\begin{lemma}
\label{lem:HAtk}
Let $A =[\mathbf{0},\BFx)$ with $\BFx \in [0,1]^s$, and let $ t_{\BFk}$ {be the coefficients for which % for which  
\eqref{eq:VAdecompElemk}  holds}. 
Then {for a scrambled $(0,m,s)-$net $\tilde{P}_n$}
\begin{equation}
\label{eq:HAPUV}
H_n(A) = \sum_{\BFk\ge \mathbf{0}} t_{\BFk} b^{2|\BFk|} 
%P((\BFU,\BFV) \in \Ik \times \Ik)).
H_n(\Ik).
\end{equation}
\end{lemma}

%is that if $A,B \subseteq [0,1]^s$ are such that $V(A) = V(B)$, then $H(A) = H(B)$. \christiane{probably need to make this into a lemma or proposition} This property, along with the decomposition \eqref{eq:VAdecompElemk} yields the identity

\begin{proof}
As shown in \cite{WLD20}, we use the (constant) value $\psi_{\BFi}$ of the joint pdf of $\BFU,\BFV$ {from $\tilde{P}_n$} over $D_{\BFi}^s$ to compute $H_n(A)$ as
\begin{equation}
\label{eq:HAinsidelemma}
H_n(A) = \sum_{\BFi\ge \mathbf{0}} \psi_{\BFi}  \times V_{\BFi}(A \times A)
\end{equation}
and then use \eqref{eq:VAdecompElemk} to get 
%write
%\begin{equation}
%    \label{eq:ViAinsidelemma}
%V_{\BFi}(A \times A) = \sum_{\BFk\ge \mathbf{0}} t_{\BFk} %b^{2|\BFk|}V_{\BFi}(\Ik \times \Ik).
%\end{equation}
%Substituting \eqref{eq:ViAinsidelemma} into \eqref{eq:HAinsidelemma} we %get
%\begin{align*}
    $H_n(A) =%&= 
    \sum_{\BFi\ge \mathbf{0}} \psi_{\BFi} \sum_{\BFk\ge \mathbf{0}}t_{\BFk} b^{2|\BFk|}V_{\BFi}(\Ik \times \Ik)
    %=\sum_{\BFk\ge \mathbf{0}} t_{\BFk} b^{2|\BFk|}
    %\sum_{\BFi\ge \mathbf{0}} \psi_{\BFi}V_{\BFi}(\Ik \times \Ik)
    %\\
    %&
    =\sum_{\BFk\ge \mathbf{0}} t_{\BFk} b^{2|\BFk|}
    H_n(\Ik)$
%    P((\BFU,\BFV) \in \Ik \times \Ik),
%\end{align*}
where the order of summation can be changed thanks to Tonelli's theorem.
\end{proof}

Next, rather than computing $H_n(\Ik)$ 
%$P((\BFU,\BFV) \in \Ik \times \Ik)$ 
by writing it as an integral involving the joint pdf %$\psi(\BFu,\BFv)$ 
associated with $(\BFU,\BFV)$ (as we just did in the proof of Lemma \ref{lem:HAtk}), we instead use a conditional probability argument that allows us to directly connect this probability to the counting numbers $m_b(\BFk;P_n)$ 
\hchristiane{explain why $m_b(\BFk;P_n)$  doesn't depend on $l$}
used in \cite{WLD20}, which for $P_n$ a digital net, represents the number of points $\BFu_j \in  P_n$ satisfying $\boldsymbol{\gamma}_b^s(\BFu_l,\BFu_j) \ge \BFk$
for a given $l \neq j$. ({For an arbitrary $P_n$, this number depends on $\ell$ but for a $(0,m,s)-$net, it is invariant with $\ell$, hence we drop  the dependence on $\ell$ in our notation. Also,} since $m_b(\BFk;P_n)=m_b(\BFk;\tilde{P}_n)$, we work with the deterministic point sets when using these counting numbers.)
{This is a key step, as it allows us to write $H_n(A)$ as a linear equation instead of an integral, thereby yielding a linear programming formulation for our main result, which is to show $H_n(A) \le {\rm Vol}(A \times A).$} Specifically, we write
\begin{align}
%P((\BFU,\BFV) \in \Ik \times \Ik)
H_n(\Ik)
&= P(\BFV \in \Ik | \BFU \in \Ik) P(\BFU \in \Ik) 
= \frac{m_b(\BFk;P_n)}{n-1} b^{-|\BFk|}. \label{eq:PUVandmk}
\end{align}
If $P_n$ is a $(0,m,s)-$net in base $b$, then  $m_b(\BFk;P_n) = \max(b^{m-|\BFk|}-1,0)$ \cite{WLD20}. %Furthermore, 
These counting numbers are also closely connected to the key quantities 
\cite{WLD20}
%$C_b(\BFk;P_n)$ from \cite{WLD20} given by
\[
C_b(\BFk;P_n) = \frac{b^{|\BFk|}m_b(\BFk;P_n)}{n-1}.
\]
%Note that these counting numbers are invariant under scrambling, which is why we can define them over the deterministic point sets $P_n$ being scrambled.

Combining \eqref{eq:HAPUV} and \eqref{eq:PUVandmk}, we get
that for $\tilde{P}_n$ a scrambled $(0,m,s)-$net,
\begin{equation}
    \label{eq:HAfromWLD20}
H_n(A) = \sum_{\BFk\ge \mathbf{0}} t_{\BFk} b^{|\BFk|} \frac{m_b(\BFk;P_n)}{n-1}
= \sum_{\BFk\ge \mathbf{0}} t_{\BFk} C_b(\BFk;P_n) 
\le \sum_{\BFk\ge \mathbf{0}} t_{\BFk}  = {\rm Vol}(A \times A),
\end{equation}
since $C_b(\BFk;P_n)  \le 1$ when $P_n$ is a  $(0,m,s)-$net \cite{WLD20}. % (also holds for the first $n$ points of a $(0,s)-$sequence.)

%%%%%%%%%%%%%%%%%%%%%%%%%%%%%%%%%%%%%
{\section{Decomposing unanchored intervals}}
\label{sec:decomp}

%[TODO: explain the approach; might need to move some proofs in an appendix]
\hchristiane{consider changing section title to ''decomposing intervals''; adapt first few sentences to explain we really just need to know how to do this in one dimension}

We now consider the case where $A$ is an unanchored box of the form $A = \prod_{j=1}^s [a_j,A_j)$, with $0 \le a_j < A_j \le 1, j=1,\ldots,s$. In Section \ref{sec:negdepunanchoredbox}, we will prove {in Theorem \ref{thm:boundHA}} that for a scrambled $(0,m,s)-$net, we still have $H_n(A) \le {\rm Vol}(A\times A)$ in this case, {which is the main result of this paper}. The proof of this result is much more difficult than in the anchored case because when $A$ is not anchored at the origin, we cannot always find a conical decomposition of products of elementary intervals as in \eqref{eq:VAdecompElemk} that has the same volume vector as $A\times A$.

Before going further, we note that it is sufficient to focus on the decomposition of one-dimensional intervals $A$ since a box is just a product of intervals. Hence for the rest of this section, we assume $s=1$.

The reason why the decomposition \eqref{eq:VAdecompElemk} cannot be used for unanchored intervals is that it may produce 
%
%For example\hjaspar{It may be better to say for example rather than to be mroe precise}, for some unanchored $A$,  the decomposition used in the anchored case leads to some 
coefficients $t_{\BFk}$ that are negative, which makes the inequality in \eqref{eq:HAfromWLD20} not necessarily true. In turn, this happens because the key property $bV_{i}(A\times A) \ge V_{i-1}(A\times A)$ {that holds for an anchored interval $A$} and that is used to show that $t_{\BFk} \ge 0$ in \cite{WLD20} is not always satisfied when $A$ is an unanchored  {interval. In this case,} 
%Before going further, we  Indeed for an unanchored one-dimensional interval $A$},  
the volume vector corresponding to $A \times A$ may be such that  $V_0(A \times A)>0$, $V_1(A \times A) = \ldots = V_{r-1}(A \times A)=0$, $V_i(A \times A) >0$ for $i \ge r$.
Because of this, we can see from \eqref{eq:deftk} that some $t_{\BFk}$ may be negative. %To address this problem, we introduce  families of vectors ${\cal W}^{(d)}$, which we now define.

To get a decomposition with non-negative coefficients, 
we introduce a family of regions of the form $Y \times Y$ where $Y$ is not an elementary interval
anchored at the origin.
More precisely, for $d,k$ non-negative integers, we define what we call an {\em elementary unanchored $(d,k)-$interval}
\[
Y_k^{(d)} := \left[ \frac{1}{b^{d+1}} - \frac{1}{b^{2+k+d}}, \frac{1}{b^{d+1}} + \frac{1}{b^{2+k+d}} \right).
\]

As a first step, in the following lemma we establish some key properties for the volume vector corresponding to an unanchored interval $A$. It is the counterpart to the property that $bV_i(A \times A) \ge V_{i-1}(A\times A)$ for anchored boxes, and shows that the $V_i(A\times A)$'s do not decrease too quickly with $i$ in the unanchored case, which is essential to prove the decomposition given in Theorem \ref{thm:decompA}. The proof of this lemma is in the appendix.
{Note that this lemma applies to half-open intervals strictly contained in $[0,1)$; the interval $[0,1)$ can be handled using the decomposition from \cite{WLD20}, which was described in the previous section.}

%%%%%%%%%%%%%%%%%%%%%%%%%%%%%%%%%%%%%%%%

\hchristiane{we need to be more precise about which $A$ does this apply to: e.g., add that $A$ must be half-open; might need to explain that the lemma is for non-anchored intervals which is why it excludes $[0,1)$ or even $hb^{-k},(h+1)b^{-k})$, as mentioned by Ref 2}

\begin{lemma}
\label{lem:vol}
Let {$A \subset [0,1)$ be a half-open interval} and let $r \ge 1$ be the smallest integer such that we can write $A = [hb^{-r+1}+gb^{-r}-z,hb^{-r+1}+Gb^{-r}+Z)$ \hjaspar{If $A=[x,y)$ then do we not have $\gamma_b(x,y)=r-1$?} \hchristiane{Only if $g \neq G$} with $0 \le h < b$, $1 \le g \le G \le b-1$ and $z,Z \in [0,b^{-r})$. %= [a_1,A_1)$ with $0 < a_1 < A_1 \le 1$. Let $r \ge 1$ be the smallest integer such that we have $a_1 = g/b^r-z$ and $A_1 = G/b^r+Z$, for some $z,Z \in (0,1/b^r)$. 
%Let $q$ be the smallest integer such that $\max(z,Z)>b^{-q}$. (Note that $q>r$.) Let $p = 0$ if $G>g$, otherwise define $p=q-r-1$ (when $g=G$). 
Then  $V(A \times A)$ is such that:
%the corresponding volume vector $V(A,A) \in \ell^1(\N_0)$ is such that:
\begin{enumerate}[label=\roman*)]
\item  $V_i(A \times A) =0$ for $i=0,\ldots,r-2$;
\item $bV_{i+1}(A \times A) \ge V_{i}(A \times A)$ for all  $i \ge r$;
\item $V_{r-1}(A \times A) - \frac{b(b-2)}{b-1}V_r(A \times A) \le \tilde{V}_r(A \times A)$, where $\tilde{V}_r(A \times A) = \sum_{i=r}^{\infty} V_i(A \times A)$.
%[NEED TO CHANGE THIS TO $b(b-2)/(b-1)$ instead of $b^2/(b-1)$: PROOF DONE ON PAPER]
%\item  $V_{r-1} > 0$; 
%\item $V_j$=0 for $j=r,\ldots,r+p-1$ 
\end{enumerate}
%where .
\end{lemma}

The next result establishes that any unanchored interval $A$ in $[0,1)$ has a volume vector $V(A \times A)$ that
can be decomposed into a conical combination of
volume vectors of elementary unanchored $(d,k)$-intervals $Y_k^{(d)}$ and elementary (anchored) $k-$intervals $1_k$.
Its proof is in the appendix.
%, the corresponding volume vector can be written as a combination of the vectors in ${\cal W}^{(r-1)}$ and ${\cal Z}$, where $r$ is defined as in Lemma \ref{lem:vol}.
%[TODO: we could state this result without making use of the volume vectors of the $Y \times Y$ regions, and just say that $A$ has the same volume vector as a combination of appropriate basis intervals. Then the main part of the paper would completely avoid having to describe the volume vectors of the $Y \times Y$ or of products of elementary intervals.]
%%%%%%%%%%%%%%%%%%%%%%%%%%%%%%%%%%%%%%%%%%%%%%%%%

\hchristiane{I think we left $d$ and it should be $r-1$, as pointed out by Ref 2: check this. Also, ref 2 points out thm would be true without the 4, so we need to explain more about the scaling, i.e., that the $\alpha_k$ are normalizing constants, and the 4 comes from cartesian square that was cut in half or something like that}
\begin{theorem}
\label{thm:decompA}
Let $A \subseteq [0,1)$ be a {half-open} interval. For $A \neq [0,1)$, let $r \ge 1$ be the smallest positive integer such that we can write $A = [hb^{-r+1}+gb^{-r}-z,hb^{-r+1}+Gb^{-r}+Z)$ with $0 \le h < b$, $1 \le g \le G \le b-1$, and $z,Z \in [0,b^{-r})$. For $A = [0,1)$, let $r=1$.
%Consider $A = [a_1,A_1)$ with $0 < a_1 < A_1 \le 1$ and $r \ge 1$ be the smallest integer such that we have $a_1 = g/b^r-z$ and $A_1 = G/b^r+Z$, for some $z,Z \in (0,1/b^r)$. 
%Let $V(A \times A)$ be the volume vector associated to $A \times A$. 
Then there exists non-negative coefficients $(\alpha_k)_{k \ge 0}$ and $(\tau_k)_{k \ge 0}$ such that {${\rm Vol}(A\times A)=\sum_{k\geq 0}(\alpha_k+\tau_k)$ and }
\[
V(A \times A) = \sum_{k = 0}^{\infty} \alpha_k \frac{b^{2(k+{r+1})}}{4} V(Y_k^{(r-1)} \times Y_k^{(r-1)}) + \sum_{k=0}^{\infty} \tau_k b^{2k} V(1_k \times 1_k).
%V = \sum_{k = 0}^{\infty} \alpha_k W_k^{(r-1)} + \sum_{k=0}^{\infty} \tau_k Z_k.
\]
\end{theorem}

%%%%%%%%%%%%%%%%%%%%%%%%%%%%%%%%%%%%%%%%%%%%%%%%
\section{Pairwise sampling dependence of scrambled $(0,m,s)-$nets on unanchored boxes}
\label{sec:negdepunanchoredbox}

{This section contains our main result,} which is 
%The main result of this section is 
that  scrambled $(0,m,s)-$nets have a negative pairwise sampling dependence index. That is, for this construction, $H_n(A) \le {\rm Vol}(A \times A)$ for any unanchored box $A$. To prove this result, we must first provide a decomposition for $H_n(A)$ 
that makes use of elementary intervals and elementary  unanchored $(d,k)-$intervals.
To do so, we use the 
decomposition of an unanchored interval given in Theorem \ref{thm:decompA}. First, we introduce some notation to denote regions in $[0,1)^{2s}$ that will be used repeatedly in this section, starting with those
%
%We start with the type of regions 
we get from the decomposition proved in Theorem \ref{thm:decompA}:
\begin{align}\label{eq:ExtreemPoints}
D(\BFk,\BFd,J) :=& 
%\prod_{j \in J} \left[\left. \frac{1}{b^{d_j+1}} - %\frac{1}{b^{d_j+k_j+2}},
%\frac{1}{b^{d_j+1}} + \frac{1}{b^{d_j+k_j+2}}
%\right.\right) \times \left[\left. \frac{1}{b^{d_j+1}} - %\frac{1}{b^{d_j+k_j+2}},
%\frac{1}{b^{d_j+1}} + \frac{1}{b^{d_j+k_j+2}}
%\right.\right)\\
\prod_{j \in J} Y_{k_j}^{(d_j)} \times Y_{k_j}^{(d_j)}
\prod_{j \in J^c} 1_{k_j} \times 1_{k_j},
%\left[\left. 0, \frac{1}{b^{k_j}} \right. \right) \times \left[\left. 0, \frac{1}{b^{k_j}} \right. \right),
\end{align}
where $J \subseteq \{1,\ldots,s\}$.
The interval $Y_{k_j}^{(d_j)}$ is decomposed further using
\[
Y_{k_j,1}^{(d_j)} := \left[\left. \frac{1}{b^{d_j+1}}-\frac{1}{b^{2+k_j+2_j}},
\frac{1}{b^{d_j+1}}\right) \right.,
\qquad \mbox{and }
Y_{k_j,2}^{(d_j)} := \left[\left. \frac{1}{b^{d_j+1}},\frac{1}{b^{d_j+1}}+
\frac{1}{b^{2+k_j+2_j}}
 \right) \right. .
\]
We also make use of the following sub-regions, where $I,K \subseteq J$:
\begin{align*}
E(\BFk,\BFd,J,I) &:= \prod_{j \in I}
1_{k_j+d_j+2} \times 1_{k_j+d_j+2} 
\prod_{j \in J^c} 1_{k_j} \times 1_{k_j} \\
%\left[\left. 0,\frac{1}{b^{d_j+k_j+2}}
%\right.\right) \times \left[\left. %0,\frac{1}{b^{d_j+k_j+2}}
%\right.\right)  
% \prod_{j \in J^c} \left[\left. 0, \frac{1}{b^{k_j}} \right. \right) \times \left[\left. 0, \frac{1}{b^{k_j}} \right. \right) %\\
\tilde{E}(\BFk,\BFd,J,I,K) &:=  \prod_{j \in I \cap K}
Y_{k_j,2}^{(d_j)} \times Y_{k_j,2}^{(d_j)}
%\left[\frac{1}{b^{d_j+1}},\frac{1}{b^{d_j+1}}+\frac{1}{b^{d_j+k+j+2}} \right)
\prod_{j \in I \cap K^c} 
%\left[\frac{1}{b^{d_j+1}}-\frac{1}{b^{d_j+k+j+2}},\frac{1}{b^{d_j+1}} \right)
Y_{k_j,1}^{(d_j)} \times Y_{k_j,1}^{(d_j)}
\prod_{j \in J^c} 1_{k_j} \times 1_{k_j}\\
F(\BFk,\BFd,J,I) &:= E(\BFk,\BFd,J,I) \times
\prod_{j \in J \cap I^c} 
Y_{k_j,1}^{(d_j)} \times Y_{k_j,2}^{(d_j)}\\
%\end{align*}
%\begin{align*}
F(\BFk,\BFd,J,I,K) &:= \tilde{E}(\BFk,\BFd,J,I,K)
\prod_{j \in J \cap I^c \cap K} 
Y_{k_j,1}^{(d_j)} \times Y_{k_j,2}^{(d_j)} 
\prod_{j \in J \cap I^c \cap K^c} 
Y_{k_j,2}^{(d_j)} \times Y_{k_j,1}^{(d_j)}.
\end{align*}

The region $F(\BFk,\BFd,J,I)$ in which a pair of points $(\BFU,\BFV)$ lies will sometimes be written as the product of the two regions obtained by projecting it over the coordinates of $\BFU$ and then $\BFV$, using the notation
$
F(\BFk,\BFd,J,I) = F_1(\BFk,\BFd,J,I) \times F_2(\BFk,\BFd,J,I).
$
That is,
$
F_i(\BFk,\BFd,J,I)  := \prod_{j \in I} 1_{k_j+d_j+2} \prod_{j \in J^c} 1_{k_j}
\prod_{j \in J \cap I^c} Y_{k_j,i}^{(d_j)}$ for $i=1,2.$
%and
%$F_2(\BFk,\BFd,J,I)  := \prod_{j \in I} 1_{k_j+d_j+2} %\prod_{j \in J^c} 1_{k_j}
%\prod_{j \in J \cap I^c} Y_{k_j,2}^{(d_j)}.
%$

Next, we define the counting numbers
$m_b(\BFk,\BFd,c,J,I;P_n)$. {The parameter $c \ge 0$ is used to specify the number of initial common digits over the subset $I$.} %As suggested from our definition of the regions $Y_k^{(d)}$ it will be set to 2 for our analysis.

\hchristiane{say something about what $c$ represents}
\begin{definition} For $P_n$ a digital net, let $m_b(\BFk,\BFd,c,J,I;P_n)$ be the number of points $\BFu_{\ell}$, for a given point $\BFu_i \in P_n$, which are different from $\BFu_i$ and satisfy: % the following properties:
\vspace*{-0.2cm}
\begin{align*}
 \gamma_b(u_{i,j},u_{\ell ,j}) &\ge k_j+d_j+c
    \mbox{ if } j \in I; \\
    \gamma_b(u_{i,j},u_{\ell,j}) &\ge k_j
    \mbox{ if } j \in J^c; \\
    \gamma_b(u_{i,j},u_{\ell,j}) &= d_j 
    \mbox{ if } j \in J\cap I^c. 
\end{align*}
\end{definition}

The properties stated in the next lemma involve the above regions and will be useful to prove Theorems \ref{thm:decompHA} and \ref{thm:boundHA}. Its proof is in the appendix.

\begin{lemma}
\label{lem:volcondprob}
Let $|\BFk|_J$ denote the sum $\sum_{j\in J} k_j$ with also $|\BFk+2|_J = \sum_{j\in J} (k_j+2) = |\BFk|_J + 2|J|$. Then:
\begin{enumerate}
\item ${\rm Vol}(D(\BFk,\BFd,J))=2^{2|J|}b^{-2(|\BFk|+|\BFd+2|_{J})}$.
\item ${\rm Vol}(F_1(\BFk,\BFd,J,I))=b^{-(|\BFk|+|\BFd+2|_{J})}$.
\item %It holds that
$
P(\BFV \in F_2(\BFk,\BFd,J,I)|\BFU \in F_1(\BFk,\BFd,J,I))
= \frac{m_b(\BFk,\BFd,2,J,I;P_n)}{n-1} \frac{(b-1)^{|I|-|J|}}{b^{|\BFk|_{J \cap I^c}+|J|-|I|}}.
$
\item We have 
%For $K \subseteq J$, let
%\[
%\tilde{E}(\BFk,\BFd,J,I,K) =  \prod_{j \in I \cap K}
%Y_{k_j,2}^{(d_j)} \times Y_{k_j,2}^{(d_j)}
%\left[\frac{1}{b^{d_j+1}},\frac{1}{b^{d_j+1}}+\frac{1}{b^{d_j+k+j+2}} \right)
%\prod_{j \in I \cap K^c} 
%\left[\frac{1}{b^{d_j+1}}-\frac{1}{b^{d_j+k+j+2}},\frac{1}{b^{d_j+1}} \right)
%Y_{k_j,1}^{(d_j)} \times Y_{k_j,1}^{(d_j)}
%\prod_{j \in J^c} 1_{k_j} \times 1_{k_j}
%\]
%and 
%
%Then 
$D(\BFk,\BFd,J) = \cup_{K,I \subseteq J} F(\BFk,\BFd,J,I,K)$
and
$P((\BFU,\BFV) \in F(\BFk,\BFd,J,I)) = 
P((\BFU,\BFV) \in F(\BFk,\BFd,J,I,K))$
for all $K,I \subseteq J$.
\end{enumerate}
\end{lemma}

%Next, we show that scrambled $(0,m,s)$-nets in base $b$ are negatively independent on intervals of the form
%\[
%\prod_{j=1}^s A_{k_j}^{(d_j)} = \left[\frac{1}{b^{d_j+1}}-\frac{1}{b^{2+d_j+k_j}},\frac{1}{b^{d_j+1}}+\frac{1}{b^{2+d_j+k_j}}\right].
%\] 

\hchristiane{There seems to be no reason why this result is stated only for scrambled $(0,m,s)$-nets: I will change this}
The next result provides us with a key decomposition for $H_n(A)$.

\begin{theorem}
\label{thm:decompHA}
Let $A = \prod_{j=1}^s [a_j,A_j)$  
be an unanchored box, %in $[0,1)^s$  of the form
where $0 \le a_j < A_j \le 1$, $j=1,\ldots,s$. 
%%Sept 25 
%%Let $\psi(\BFu,\BFv)$ be the joint pdf of a scrambled $(0,m,s)$-net in base $b$. 
For $J \subseteq\{1,\ldots,s\}$ and {$\BFk \in \N_0^s$} , let 
 $\tilde{\alpha}_{\BFk,J} = \prod_{j \in J} \alpha_{k_j}^{(j)}
   \prod_{j \in J^c} \tau_{k_j}^{(j)},$
where the %coefficients 
$\alpha_{k_j}^{(j)}$ and $\tau_{k_j}^{(j)}$ come from the decomposition given in Theorem \ref{thm:decompA} applied to the interval $[a_j,A_j)$, $j=1,\ldots,s$. In particular, this means
$\tilde{\alpha}_{\BFk,J}\ge 0$ and
\begin{equation}
    \label{eq:VolInTermsCoeff}
 \sum_{\BFk \ge \mathbf{0}} \sum_{J \subseteq \{1,\ldots,s\}}  \tilde{\alpha}_{\BFk,J} = {\rm Vol}(A \times A).   
\end{equation}

Let $(\BFU,\BFV)$ be a randomly chosen pair of points from a point set $\tilde{P}_n$.
Then
\[
H_n(A) %= \int_{A \times A} \psi(\BFu,\BFv) d\BFu d\BFv
= \sum_{\BFk \ge \mathbf{0}} \sum_{J \subseteq \{1,\ldots,s\}}  \frac{\tilde{\alpha}_{\BFk,J}}{{\rm Vol}(D(\BFk,\BFd,J))}
P((\BFU,\BFV) \in D(\BFk,\BFd,J)).
%%Sept 25 \tilde{C}_{\BFk,J},
%\le  \sum_{\BFk} \sum_{J}  \tilde{\alpha}_{\BFk,J} = {\rm Vol}(A \times A).
\]
%where 
%\begin{align*}
%   \tilde{\alpha}_{\BFk,J} &= \prod_{j \in J} %\alpha_{k_j}^{(j)}
%   \prod_{j \in J^c} \tau_{k_j}^{(j)} \\
%   \tilde{C}_{\BFk,J} &= \sum_{\BFi \ge \mathbf{0}} \psi_{\BFi}\prod_{j \in J} W_{i_j,k_j}^{(d_j)} \prod_{j \notin J} Z_{i_j,k_j}
%\end{align*}
%are both scalars \hjaspar{Consider adding this, or something like it, for clarity.} and where $W_{\BFk}^{(\BFd)}$ is defined in \eqref{eq:Wikd} and $Z_{\BFk}$ is defined in \eqref{eq:zik}, while the coefficients $\alpha_{k_j}^{(j)}$ and $\tau_{k_j}^{(j)}$ come from the decomposition studied in Theorem \ref{thm:decompA} applied to the interval $[a_j,A_j]$, $j=1,\ldots,s$. Further, it can be shown that
\end{theorem}

\begin{proof}
Using Theorem \ref{thm:decompA}  we can write
\hchristiane{Ref 2 saying that $V$ should only have one argument and that $D(\BFk,\BFd,J)$ cannot be an argument of $V$: check all that}
\begin{align*}
V(A \times A) &= \prod_{j=1}^s \left(\sum_{k_j \ge 0} \alpha_{k_j}^{(j)} \frac{b^{2(k_j+d_j+2)}}{4} V(Y_{k_j}^{(d_j)} \times Y_{k_j}^{(d_j)})+ \sum_{k_j \ge 0} \tau_{k_j}^{(j)}b^{2k_j} 
V(1_{k_j} \times 1_{k_j})\right)\\
&=\sum_{\BFk \ge \mathbf{0}} \sum_{J \subseteq \{1,\ldots,s\}}
\left(\prod_{j \in J} \alpha_{k_j}^{(j)} \frac{b^{2(k_j+d_j+2)}}{4}  \prod_{j \in J^c} b^{2k_j}\tau_{k_j}^{(j)}\right) V(D(\BFk,\BFd,J)) \\
&=\sum_{\BFk \ge \mathbf{0}} \sum_{J \subseteq \{1,\ldots,s\}} \tilde{\alpha}_{\BFk,J} 2^{-2|J|}b^{2(|\BFk|+|\BFd+2|_J)} V(D(\BFk,\BFd,J))\\
&=\sum_{\BFk \ge \mathbf{0}} \sum_{J \subseteq \{1,\ldots,s\}} \frac{\tilde{\alpha}_{\BFk,J}}{{\rm Vol}(D(\BFk,\BFd,J))}  V(D(\BFk,\BFd,J)),
\end{align*}
where the last equality follows from 
Part 1 of Lemma \ref{lem:volcondprob}. Then, using the same kind of reasoning as in Lemma \ref{lem:HAtk}, we get
\begin{align*}
H_n(A) &=\sum_{\BFk \ge \mathbf{0}} \sum_{J \subseteq \{1,\ldots,s\}} \frac{\tilde{\alpha}_{\BFk,J}}{{\rm Vol}(D(\BFk,\BFd,J))} 
P((\BFU,\BFV) \in D(\BFk,\BFd,J)). \,\,\,\,\,\,\hspace*{1cm} \square
\end{align*}
%as required. 
\end{proof}

It is clear from Theorem \ref{thm:decompHA} that in order to prove that $H_n(A) \le {\rm Vol}(A \times A)$, it is sufficient to prove that 
$
P((\BFU,\BFV) \in D(\BFk,\BFd,J)) \le {\rm Vol}(D(\BFk,\BFd,J))
$
for all $\BFk,\BFd,J.$
{That is, the regions $D(\BFk,\BFd,J)$ correspond to the extreme points in the linear programming formulation of our problem.}

A key quantity to analyze this probability is the following weighted sum of counting numbers for $P_n$, where $J \subseteq \{1,\ldots,s\}$ and $I^* := I \cup J^c$:
\begin{equation}
\label{eq:weightedCountNbers}
\tilde{m}_b(\BFk,\BFd,J;P_n) := 
\frac{1}{2^{|J|}}\!\!\! \sum_{I \subseteq J}b^{|\BFk|_{I^*}+|\BFd|_J+|J|+|I|}
(b-1)^{|I|-|J|}
\frac{m_b(\BFk,\BFd,2,J,I;P_n)}{n-1}.
\end{equation}

\begin{theorem}
\label{thm:boundHA}
Let $A$ be an unanchored box in $[0,1)^s$. Let  $H_n(A)$ and  $\tilde{\alpha}_{\BFk,J}$ %, and $\tilde{C}_{\BFk,J}$ 
be defined as in Theorem \ref{thm:decompHA}. 
Let $P_n$ have counting numbers $m_b(\BFk,\BFd,2,J,I;P_n)$ such that
\begin{equation}
\label{ConditionWeightedCountNbers}
\tilde{m}_b(\BFk,\BFd,J;P_n)
\le 1,
\end{equation}
where $\tilde{m}_b(\BFk,\BFd,J;P_n)$ is defined in 
\eqref{eq:weightedCountNbers}.
Then {the scrambled point set $\tilde{P}_n$ is such that}
\[
H_n(A)%= \sum_{\BFk} \sum_{J}  %\tilde{\alpha}_{\BFk,J}\tilde{C}_{\BFk,J}
\le  \sum_{\BFk \ge \mathbf{0}} \sum_{J \subseteq \{1,\ldots,s\}}  \tilde{\alpha}_{\BFk,J} = {\rm Vol}(A \times A).
\]
\end{theorem}

\begin{proof}
As mentioned earlier, based on Theorem \ref{thm:decompHA}, it suffices to show that
%%Sept 25 $\tilde{C}_{\BFk,J}\le 1$ 
$
P((\BFU,\BFV) \in D(\BFk,\BFd,J)) \le {\rm Vol}(D(\BFk,\BFd,J))
$
for all 5-tuples $(m,s,\BFk,\BFd,J)$, where $m \ge 1, s \ge 1, \BFk \ge \mathbf{0},\BFd \ge \mathbf{0},J \subseteq \{1,\ldots,s\}$. Indeed, if this holds, then from Theorem \ref{thm:decompHA} and using the fact that  $\tilde{\alpha}_{\BFk,J} \ge 0$, we can derive the inequality
\[
H_n(A)= \sum_{\BFk} \sum_{J}  \tilde{\alpha}_{\BFk,J}
\frac{P((\BFU,\BFV) \in D(\BFk,\BFd,J))}{{\rm Vol}(D(\BFk,\BFd,J))} 
\le \sum_{\BFk} \sum_{J}  \tilde{\alpha}_{\BFk,J} = {\rm Vol}(A \times A),
\]
where the last equality is obtained from  \eqref{eq:VolInTermsCoeff}, also proved in Theorem \ref{thm:decompHA}.

To analyze the probability $P((\BFU,\BFV) \in D(\BFk,\BFd,J))$, we use the decomposition of 
$D(\BFk,\BFd,J)$ into the sub-regions $F(\BFk,\BFd,J,I)$
outlined in Part 4 of Lemma \ref{lem:volcondprob}:
\hchristiane{will need to explain where the $2^{|J|}$ comes from. Essentially the number of copies of the $F()$ regions we need to get the $D()$ region.}
\begin{align}
&\frac{P((\BFU,\BFV) \in D(\BFk,\BFd,J))}{{\rm Vol}(D(\BFk,\BFd,J))}
= \sum_{I \subseteq J} 2^{|J|}
\frac{P((\BFU,\BFV) \in F(\BFk,\BFd,J,I))}{{\rm Vol}(D(\BFk,\BFd,J))} \notag\\
&=\frac{1}{{\rm Vol}(D(\BFk,\BFd,J))} 
\sum_{I \subseteq J} 2^{|J|}
P(\BFU \in F_1(\BFk,\BFd,J,I))
P(\BFV \in F_2(\BFk,\BFd,J,I)|\BFU \in F_1(\BFk,\BFd,J,I))\notag \\
&=\sum_{I \subseteq J} 2^{|J|} \frac{{\rm Vol}(F_1(\BFk,\BFd,J,I))}{{\rm Vol}(D(\BFk,\BFd,J))}
\frac{m_b(\BFk,\BFd,2,J,I;P_n)}{n-1} \frac{(b-1)^{|I|-|J|}}{b^{|\BFk+1|_{J \cap I^c}}} \notag \\
&=\sum_{I \subseteq J} 2^{|J|} \frac{b^{-(|\BFk|+|\BFd+2|_J)}}{2^{2|J|}b^{-2(|\BFk|+|\BFd+2|_J)}}\frac{m_b(\BFk,\BFd,2,J,I;P_n)}{n-1} \frac{(b-1)^{|I|-|J|}}{b^{|\BFk+1|_{J \cap I^c}}} \notag  
\end{align}
\begin{align}
&=\sum_{I \subseteq J} \frac{1}{2^{|J|}}
b^{|\BFk|_{I^*}+|\BFd|_J+|J|+|I|}(b-1)^{|I|-|J|} \frac{m_b(\BFk,\BFd,2,J,I;P_n)}{n-1} 
=\tilde{m}_b(\BFk,\BFd,J;P_n)\le 1,
%\label{eq:TheProb}
\notag
\end{align}
where the first equality comes from Lemma \ref{lem:volcondprob} (Part 4), the third from Lemma \ref{lem:volcondprob} (Part 3), the fourth from Lemma \ref{lem:volcondprob} (Parts 1, 2), and the last inequality follows from  \eqref{ConditionWeightedCountNbers}.
\end{proof}

{To get to our ultimate goal---which is captured in Theorem \ref{thm:scr0msnegdepunanchored} and is to prove that   $H_n(A) \le {\rm Vol}(A \times A)$ for an unanchored box $A$ for a scrambled $(0,m,s)-$net---thanks to Theorem \ref{thm:boundHA} all we need to do is to show that}
%Based on this result, what remains to be shown \christiane{be more clear that the goal is Thm 4, which is many pages away: might be solved if we're more clear in the intro about the main result and have a sketch of the proof} is that 
the condition \eqref{ConditionWeightedCountNbers}
indeed holds for a  $(0,m,s)-$net. The rest of this section is devoted to this {(cumbersome)} task.

First we write $\tilde{m}_b(\BFk,\BFd,J;P_n) = \sum_I \psi_m(\BFk,\BFd,J,I)/2^{|I|}$, where
%$\tilde{m}_b(\BFk,\BFd,J;P_n)$ as an average of terms of the form
\begin{equation}
    \label{eq:theterms}
\psi_m(\BFk,\BFd,J,I) := b^{|\BFk|_{I^*}+|\BFd|_J+|J|+|I|}(b-1)^{|I|-|J|} \frac{m_b(\BFk,\BFd,2,J,I;P_n)}{n-1}.
\end{equation}

The difficulty that arises when trying to bound the sum \eqref{eq:weightedCountNbers} by 1 is that some of the terms \eqref{eq:theterms} can be larger than 1 for certain combinations of $m,\BFk,\BFd$, and $J$. Hence we need to show that the smaller terms compensate for those larger than 1 so that overall, the average of these terms is indeed bounded by 1.

Now, we will not work directly with the counting numbers
$m_b(\BFk,\BFd,2,J,I;P_n)$ and will instead bound them, which in turn will yield a bound on $\psi_m(\BFk,\BFd,J,I)$ via \eqref{eq:theterms} and thus a bound on $\tilde{m}_b(\BFk,\BFd,J;P_n)$.
To show this bound is no larger than 1.
%Since $\tilde{m}_b(\BFk,\BFd,J;P_n) = \sum_I \psi_m(\BFk,\BFd,J,I)/2^{|I|}$, this means we will then get a bound on $\tilde{m}_b(\BFk,\BFd,J;P_n)$, which we will then show is no larger than 1. 
%
%This will be done by 
we break the problem in different cases, depending on the relative magnitude of $m$ vs.\  $|\BFk|,|\BFd|_J$ and $|J|$, with resulting bounds shown in Propositions \ref{prop:mkdjeasycase}, \ref{prop:mdkjsmallm} and \ref{prop:mkdjmediumm}. 

The bounds on the $\psi_m(\BFk,\BFd,I,J)$ terms will make use of the following functions.

\begin{definition}
Let $\ell,j,i$ be non-negative integers with $j>i$.
We define
\begin{equation}
    \label{eq:defhijl}
h_{j,i}(\ell) = \frac{b^{j+i-\ell}}{(b-1)^{j-i}} \binom{j-i-1}{\ell-2i}, \qquad 2i+1 < \ell < j+i, \qquad 0 \le i < j,
\end{equation}
%and
\begin{equation}
\label{eq:Defgijl}
\mbox{ and } g_{j,i}(\ell) =
\begin{cases}
1  & \mbox{ if  $\ell \ge i+j$ or, if  $\ell > 2i$ and $\ell$ is even}\\
 1+ h_{j,i}(\ell) &  \mbox{ if $2i+1 < \ell < j+i$ and $\ell$ is odd}\\
\left(\frac{b}{b-1}\right)^{j-i-1} & \mbox{ if $\ell=2i+1$}\\
0 & \mbox{ if } \ell \le 2i.
\end{cases}
\end{equation}
\end{definition}

In some cases, the following bound on $g_{i,j}(\ell)$ will be enough for our purpose. (Both Lemmas \ref{gjiloosebound}
and \ref{lem:boundPsi} are proved in the appendix.) %\ref{sec:proofs}.
\begin{lemma}
\label{gjiloosebound}
Let $j>i \ge 0$. Then 
\[
g_{i,j}(\ell) \le \left(\frac{b}{b-1}\right)^{i+j-\ell}
\mbox{ when } 2i < \ell < i+j.
\]
\end{lemma}

The next lemma gives a bound on $\psi_m(\BFk,\BFd,J,I)$
in the case of a $(0,m,s)-$net.

\begin{lemma} 
\label{lem:boundPsi}
If $P_n$ is a $(0,m,s)-$net, then 
for $I \subset J$,  $\psi_m(\BFk,\BFd,J,I)$ satisfies 
\begin{equation}
    \label{eq:BoundPsiByg}
\psi_m(\BFk,\BFd,J,I)  \le \frac{b^{m}}{b^m-1}   g_{|J|,|I|}(m-|\BFk|_{I^*}-|\BFd|_J).
\end{equation}
Moreover, when   $2|I| < m-|\BFk|_{I^*}-|\BFd|_J < |J|+|I|$, then
\begin{equation}
\label{eq:g2bound}
\psi_m(\BFk,\BFd,J,I)  \le \frac{b^m}{b^m-1} \left(\frac{b-1}{b}\right)^{m-|\BFk|_{I^*}-|\BFd|_J-|I|-|J|}.
\end{equation}
\end{lemma}

%%%%%%%%%%%%%%%%%%%%%%%%%%%%%%%%

Having found  a bound %$g_{|J|,|I|}(m-|\BFk|_{I^*}-|\BFd|_J)$ 
for $\psi_m(\BFk,\BFd,J,I) $  for the possible ranges of values for $m$, we can now return to the task of bounding the weighted sum $\tilde{m}_b(\BFk,\BFd,J;P_n)$.
We start with the easiest case.

\begin{proposition}
\label{prop:mkdjeasycase}
Let $P_n$ be a $(0,m,s)-$net.
If $J \neq \emptyset$ and $m \ge |\BFk|+|\BFd|_J+2|J| $ then
$\tilde{m}_b(\BFk,\BFd,J;P_n) \le 1$.
\end{proposition}

\begin{proof}

In this case, $m-|\BFk|_{I^*}-|\BFd|_J \ge |I|+|J|$ for all $I$ (for a given $J$) and therefore 
\begin{align*}
\tilde{m}_b(\BFk,\BFd,J;P_n) &= \frac{1}{2^{|J|}}
\sum_{I \subseteq J} \psi_m(\BFk,\BFd,J,I)
\le \frac{1}{2^{|J|}} \frac{b^m}{b^m-1}\sum_{I \subset J}  1  
+\frac{1}{2^{|J|}} \psi_m(\BFk,\BFd,J,J),
\end{align*}
where the inequality is derived from Lemma \ref{lem:boundPsi}. Observing that
$m_b(\BFk,\BFd,2,J,J;P_n) = m_b(\tilde{\BFk};P_n)$,
where $\tilde{k}_j := k_j$ if $j \in J^c $
and $\tilde{k}_j := k_j+d_j+2$ if $j \in J$,
we get
\begin{equation}
\label{eq:psiJJ}
\psi_m(\BFk,\BFd,J,J)
= b^{|\BFk|+|\BFd|_J+2|J|} \frac{m_b(\tilde{\BFk};P_n)}{n-1}
= b^{|\BFk|+|\BFd|_J+2|J|}\frac{b^{m-|\BFk|-|\BFd|_J-2|J|}-1}{n-1}
\end{equation}
and therefore obtain
\begin{align*}
\tilde{m}_b(\BFk,\BFd,J;P_n)  &\le \frac{1}{2^{|J|}} \frac{b^m}{b^m-1}\sum_{I \subset J}  1   
 +   \frac{1}{2^{|J|}}  b^{|\BFk|+|\BFd|_J+2|J|}  \frac{b^{m-|\BFk|-|\BFd|_J-2|J|}-1}{b^m-1}\\
&= \frac{2^{|J|}-1}{2^{|J|}} \frac{b^m}{b^m-1}
 + \frac{1}{2^{|J|}}   \frac{b^m-b^{|\BFk|+|\BFd|_J+2|J|}}{b^m-1}\\
   &=  \frac{1}{2^{|J|}(b^m-1)} \left(b^m 2^{|J|} - b^{|\BFk|+|\BFd|_J+2|J|}\right)\\
   &= \frac{b^m}{b^m-1} \left(\frac{2^{|J|}-b^{|\BFk|+|\BFd|_J+2|J|-m}}{2^{|J|}}\right).
%&=  \left(\frac{b}{2(b-1)}\right)^s \frac{1}{b^m-1} 
\end{align*}
%which gives us the first line of the bound %\eqref{eq:boundCkj}.
Therefore $\tilde{m}_b(\BFk,\BFd,J;P_n) \le 1$ if 
%\begin{align*}
$b^m 2^{|J|} - b^{|\BFk|+|\BFd|_J+2|J|} \le 2^{|J|} (b^m-1)$,
or equivalently, if %\Leftrightarrow  
$2^{|J|} \le  b^{|\BFk|+|\BFd|_J+2|J|}$,
%\end{align*}
which is true since $b \ge 2$,  and $|\BFk|+|\BFd|_J+2|J| %\ge 2|J| 
> |J|$.
\end{proof}

%\noindent {\bf Case 2:}  $m< |\BFk|+|\BFd|_J+2|J|$. 

Next, we deal with the more difficult case
$m< |\BFk|+|\BFd|_J+2|J|$, which implies that the bound given in Lemma \ref{lem:boundPsi} for $\psi_m(\BFk,\BFd,J,I)$
is sometimes larger than 1.
Note that  from \eqref{eq:psiJJ}, we see that $m_b(\tilde{\BFk};P_n)=0$
and thus 
$\psi_m(\BFk,\BFd,J,J)=0$ in this case. 
%This is easy to see from , because in this case

To handle this case, we need to analyze the function
\[
G(m,s,J,\BFk,\BFd) = \sum_{I \subset J} g_{|J|,|I|}(m-|\BFk|_{I^*}-|\BFd|_J),
\]
which we may at times write as 
\[
G(m,s,J,\BFk,\BFd) = \sum_{I \subset J:m^* > 2|I| +1}\!\! \!\!\!\!1 + \sum_{I \in {\cal M}(J)} 
\!\! h_{|J|,|I|}(m^*) + \sum_{I\subset J:m^*=2|I| +1}\!\! \left(\frac{b}{b-1}\right)^{|J|-0.5(m^*-1)-1}
\]
where $m^* := m-|\BFk|_{I^*}-|\BFd |_J$ and ${\cal M}(J) = \{I\subset J:2|I| +1 < m^* < |I|+|J|, m^* \mbox{odd}\} $. 

To show 
$\tilde{m}_b(\BFk,\BFd,J;P_n) =\sum_{I \subset J}\psi_m(\BFk,\BFd,J,I)/2^{|J|}\le 1$, 
from the bound \eqref{eq:BoundPsiByg} on $\psi_m(\BFk,\BFd,J,I)$ we see it  is sufficient to show  $G(m,s,J,\BFk,\BFd) \le 2^{|J|}\frac{b^m-1}{b^m}$ since then
\[
\tilde{m}_b(\BFk,\BFd,J;P_n) \le \frac{1}{2^{|J|}}
\frac{b^m}{b^m-1} G(m,s,J,\BFk,\BFd)
\le 1.
\]

%\begin{remark}
%\label{rem:smoins2}
%If $I$ is such that $|I| = |J|-1$ then we cannot have %$2|I| < m^* < |I|+|J|$ since in that case
%$|I|+|J| - 2|I| = |J|-|I| = 1$.
%\end{remark}

The following lemma will allow us to set $\BFd=\mathbf{0}$ when  bounding $G(m,s,J,\BFk,\BFd)$.

%\noindent{\bf Claim:} for a given $m$, $s$ and $\BFk$, the %right-hand side of  \eqref{eq:BoundCase2}  
%function $G(m,s,\BFk,\BFd)$ is maximized when $\BFd= \mathbf{0}$  [TODO: write out the proof based on shift]. 
\begin{lemma} If $g_0\ge 0$ is a constant such that $G(m,s,J,\BFk,\mathbf{0}) \le g_0$ for all $(m,s,J,\BFk)$, then  $G(m,s,J,\BFk,\BFd) 
\le g_0 $ for all $(m,s,J,\BFk,\BFd)$.
\end{lemma}

\begin{proof}
If $m < |\BFd|_J$ then $m^*< 2|I|$ for all $I \subset J$ and therefore $G(m,s,J,\BFk,\BFd)=0$. If $m \ge |\BFd|_J$ then it is easy to see that $G(m,s,J,\BFk,\BFd) = G(m-|\BFd|_J,s,J,\BFk,\mathbf{0})$, because $(m,s,J,\BFk,\BFd)$ and $(m-|\BFd|_J,s,J,\BFk,\mathbf{0})$ yield the same $m^*$ for all $I \subset J$, and $G(m,s,J,\BFk,\BFd)$ only depend on $m$, $\BFk$, and $\BFd$ through $m^*$.
\end{proof}

Based on this result, we set $\BFd=\mathbf{0}$ in what follows, and
%prove a bound $G(m,s,J,\BFk,\mathbf{0})$ that will 
consider two different sub-cases. The respective bounds on $G(m,s,J,\BFk,\mathbf{0})$  are given in Propositions \ref{prop:mdkjsmallm} and \ref{prop:mkdjmediumm}, which also establish that the condition \eqref{ConditionWeightedCountNbers}---stating that $\tilde{m}_b(\BFk,\BFd,J;P_n) \le 1$---holds for each sub-case. Before we state and prove these two propositions, we first state a technical lemma needed in the proof of Proposition \ref{prop:mdkjsmallm}, and proved in the appendix.

\begin{lemma}
\label{lem:Mbs}
For $b \ge s \ge 2$ and $\tilde{s} = \lfloor \frac{s}{2} \rfloor -1$. Then 
\[
R(b,s) := \frac{1}{2^s}\sum_{j=0}^{\tilde{s}} \binom{s}{j}  \left( \frac{b}{b-1}\right)^{s-j} \le 1.
\]
%Then $R(b,s) \le 1$.
\end{lemma}

\begin{proposition}
\label{prop:mdkjsmallm}
Let $P_n$ be a $(0,m,s)-$net. 
If $J \neq \emptyset$ and $m<|J|$ then $G(m,s,J,\BFk,\mathbf{0}) \le (b-1)/b$
for all $\BFk$ and therefore $\tilde{m}_b(\BFk,\mathbf{0},J;P_n) \le 1.$
\end{proposition}

\begin{proof}
The fact that $m<|J|$ implies %the first sum in $G(m,s,\BFk,\mathbf{0})$ is 0, since we then 
$m - |\BFk|_{I^*}<|J|+|I| $ for all $I$. Also, if  $|I| \ge 0.5(|J|-1)$, then  $m - |\BFk|_{I^*}\le 2|I|$ for all $\BFk$  and then $g_{|J|,|I|}(m-|\BFk|_{I^*})=0$. Thus    %the bound in \eqref{eq:BoundCase2} is given by
\[
G(m,s,J,\BFk,\mathbf{0}) \le  \sum_{I:|I| < 0.5(|J|-1)} g_{|J|,|I|}(m-|\BFk|_{I^*}).
\]
It turns out that in this case, the simpler but larger bound \eqref{eq:g2bound} %$g_2(m,I,\BFk,\mathbf{0} )$ 
can be used (since the only non-zero $g_{|J|,|I|}(m-|\BFk|_{I^*})$ terms are those for which $2|I| < m -|\BFk|_{I^*}< |I|+|J|$, which means \eqref{eq:g2bound} can indeed be applied), so we have
\begin{align}
%\psi \cdot W^{(\mathbf{0})}_{\BFk} 
G(m,s,J,\BFk,\mathbf{0}) &\le %\frac{1}{2^s} \frac{b^m}{b^m-1} 
%\left( 
\sum_{I:|I| < 0.5(|J|-1)} \left(\frac{b-1}{b}\right)^{m-|\BFk|_{I^*}-|J|-|I|} 
 \le %\frac{1}{2^s} \frac{b^m}{b^m-1} \left( 
 \sum_{i=0}^{\lfloor 0.5|J| \rfloor -1}   \binom{s}{i}\left(\frac{b}{b-1}\right)^{|J|-i-1}\notag\\
 &= \frac{b-1}{b}  \sum_{i=0}^{\lfloor 0.5|J| \rfloor -1}   \binom{|J|}{i}\left(\frac{b}{b-1}\right)^{|J|-i},
 \label{eq:Gsmallm}
 % \right)
\end{align}
where the second inequality comes from the fact that $m-|\BFk|_{I^*}> 2|I|$ implies $|J|+|I|+|\BFk|_{I^*}-m \le |J|-|I|-1$.

Using Lemma \ref{lem:Mbs} with $s=|J|$, we get that
the sum in \eqref{eq:Gsmallm} is bounded by 1. Hence 
\[
\tilde{m}_b(\BFk,\mathbf{0},J;P_n) \le  \frac{1}{2^{|J|}} \frac{b^m}{b^m-1}  \frac{b-1}{b} < 1, \,\, \mbox{for any }|J| \ge 2, m \ge 1.  \,\,\, \square
\]

%This yields the third case given in the bound \eqref{eq:boundCkj}.
\end{proof}

The last case we need to deal with is when $m$ is such that  $|J| \le m < |\BFk|+2|J|$. % and $m$ is odd, where $k_{(1)}$ and $k_{(2)}$ represent the smallest and second smallest value of $\BFk$, respectively. 
Let ${\cal B}$ be the set of pairs $(m,\BFk)$ satisfying this assumption. 
%This corresponds to the case where the function $G(m,s,\BFk,\mathbf{0})$ has non-zero terms in  the first and (in most cases, unless $m$ is very close to its upper bound of $|\BFk|+2s-1$) in the second sum in \eqref{eq:fctG}. 

To handle this case, we make use of the following two lemmas about  $G(m,s,\BFk,J,\mathbf{0})$. 
The first one shows that when $\BFk = \BFd = \mathbf{0}$ the maximum is reached when $m=2|J|-1$.
The second one shows it is sufficient to bound 
$G(m,s,J,\BFk,\mathbf{0})$ at $\BFk = \mathbf{0}.$

\begin{lemma}
\label{lem:step1}
If $\BFk = \BFd = \mathbf{0}$, then
$G(m,s,J,\mathbf{0},\mathbf{0}) \le G(2|J|-1,s,J,\mathbf{0},\mathbf{0})
=2^{|J|}-1$ for all
$m$ such that $(m,\mathbf{0}) \in {\cal B}$.
\end{lemma}

\begin{lemma}
\label{lem:step2}
Consider a pair $(m,\BFk)$ with possibly $\BFk \neq \mathbf{0}$. %such that 
%$G(m,\BFk) \ge G(m',\BFk')$ for all $(m',\BFk') \in {\cal B}$. 
 Then there  exists an odd integer value $\tilde{m}$ such that  $G(\tilde{m},s,J,\mathbf{0},\mathbf{0}) \ge G(m,s,J,\BFk,\mathbf{0})$.
\end{lemma}

Using these two lemmas (proved in the appendix), we get a bound on 
$G(m,s,J,\BFk,\mathbf{0})$ for this last case, which in turn allows us to show that \eqref{ConditionWeightedCountNbers} also holds then.

\begin{proposition}
\label{prop:mkdjmediumm}
Let $P_n$ be a $(0,m,s)-$net. 
Assume $m$ is such that $|J| \le m < |\BFk| + 2|J|$.
Then $G(m,s,J,\BFk,\mathbf{0}) \le 2^{|J|}-1$ and therefore $\tilde{m}_b(\BFk,\mathbf{0},J;P_n) \le 1.$
\end{proposition}

\begin{proof}
For a given $s$ and $J$, % and assuming $\BFd = \mathbf{0}$, 
we need to find a bound for %the function  
$G(m,s,J,\BFk,\mathbf{0})$
over all pairs $(m,\BFk) \in {\cal B}$, and do so by showing it is 
%We will show that  this function  is 
maximized when $\BFk= \mathbf{0}$ and $m=2|J|-1$.

First, from Lemma \ref{lem:step2}, we have that 
for a given $(m,\BFk) \in {\cal B}$, we can find  a pair in ${\cal B}$ of the form $(\tilde{m},\mathbf{0})$ such that $G(\tilde{m},s,J,\mathbf{0},\mathbf{0})\ge G(m,s,J,\BFk,\mathbf{0})$. Hence we can set $\BFk = \mathbf{0}$. Next, we use Lemma \ref{lem:step1}, which shows that for pairs in ${\cal B}$ of the form $(m,\mathbf{0})$, the function $G(m,s,J,\BFk,\mathbf{0})$ is maximized when $m = 2|J|-1$.

Putting these two lemmas together, we get that for a given $s$ and $J$, $G(m,s,J,\BFk,\mathbf{0}) \le G(\tilde{m},s,J,\mathbf{0},\mathbf{0}) \le G(2|J|-1,s,J,\mathbf{0},\mathbf{0})=2^{|J|}-1$ for all $(m,\BFk) \in {\cal B}$.
Hence
\[
\tilde{m}_b(\BFk,\BFd,J;P_n) \le \frac{1}{2^{|J|}} \frac{b^{m}}{b^{m}-1} (2^{|J|}-1)
=  \frac{b^{m}}{b^{m}-1} \frac{2^{|J|}-1}{2^{|J|}} \le 1,
\]
which holds since $2^{|J|} \le b^m$, as $b \ge 2$, and $m\ge |J|$. 
\end{proof}

%We now have all that is needed to show that the %condition
%\eqref{ConditionWeightedCountNbers} holds.
Having examined all possible cases, we can now state our main result.

\begin{theorem}
\label{thm:scr0msnegdepunanchored}
%The condition \eqref{ConditionWeightedCountNbers}
%holds for a $(0,m,s)-$net.
%Therefore 
If $\tilde{P}_n$ is a scrambled $(0,m,s)-$net in base $b$,  then $H_n(A) \le {\rm Vol}(A\times A)$ for any unanchored box $A \in {\cal A}$, and thus its pairwise sampling dependence index satisfies ${\cal E}_n(\tilde{P}_n) \le 0$.
\end{theorem}

\begin{proof}
Using Theorem \ref{thm:boundHA}, we need to show that $\tilde{m}_b(\BFk,\BFd,J;P_n) \le 1$ for all $\BFk,\BFd,J$ for a $(0,m,s)-$net $P_n$, {i.e., that condition \eqref{ConditionWeightedCountNbers}
holds for a $(0,m,s)-$net.}
First, from Proposition \ref{prop:mkdjeasycase}, if $J \neq \emptyset$ and $m \ge |\BFk|+|\BFd|_J+2|J|$ then
\[
\tilde{m}_b(\BFk,\BFd,J;P_n) 
\le
\frac{b^m}{b^m-1} \left(\frac{2^{|J|} - b^{|\BFk|+|\BFd|_J+2|J|-m}}{2^{|J|}}\right)
\le 1.
\]

Next, from Proposition \ref{prop:mdkjsmallm} we have that if $J \neq \emptyset$ and $m<|J|$ then
\[
\tilde{m}_b(\BFk,\BFd,J;P_n) 
\le
\frac{b^m}{b^m-1} \frac{1}{2^{|J|}} \frac{b-1}{b}
\le 1.
\]

Then, using Proposition \ref{prop:mkdjmediumm} we get that if $0< |J|<m<|\BFk|+|\BFd|_J+2|J|$ then
\[
\tilde{m}_b(\BFk,\BFd,J;P_n) 
\le
\frac{b^m}{b^m-1} \frac{2^{|J|}-1}{2^{|J|}}
\le 1.
\]
Finally, if $J =\emptyset$ then 
\begin{align*}
\tilde{m}_b(\BFk,\BFd,J;P_n)  &= \frac{P((\BFU,\BFV) \in D(\BFk,\BFd,J))}{{\rm Vol}(D(\BFk,\BFd,J))}
=\frac{1}{b^{-2|\BFk|}}b^{-|\BFk|}\frac{\max(b^{m-|\BFk|}-1,0)}{b^m-1} \\
&= b^{|\BFk|}\frac{\max(b^{m-|\BFk|}-1,0)}{b^m-1}
=C_b(\BFk;P_n),
\end{align*}
which was shown to be smaller or equal to 1 in \cite{WLD20} for a  $(0,m,s)-$net.
\end{proof}
% end of proof of Theorem \ref{thm:boundHA}

Using Theorem \ref{thm:scr0msnegdepunanchored}, we obtain the following result, which shows that a scrambled net integrates the indicator function $1_A$ of any unanchored box $A$ with variance no larger than the Monte Carlo estimator variance.
To our knowledge, this was not previously known. What is well known is that since $A$ is an axis-parallel box, then $1_A$ has bounded variation in the sense of Hardy and Krause and therefore a scrambled net has variance in $O(n^{-2} (\log n)^s)$ \cite{vOWE97c}.

\begin{proposition}
Let $A$ be an unanchored box in $[0,1)^s$. Let $\hat{\mu}_{n,A}$ be the estimator for $\mu_A = {\rm E}(1_A) = {\rm Vol}(A)$ based on a scrambled $(0,m,s)-$net in base $b$ with $n=b^m$.
Then ${\rm Var}(\hat{\mu}_{n,A}) \le \mu_A(1-\mu_A)/n$.
\end{proposition}

\begin{proof}
The result follows from the fact that ${\rm Var}(\hat{\mu}_{n,A}) = \mu_A(1-\mu_A)/n +  (H_n(A)-\mu^2_A) (n-1)/n$, and then applying Theorem \ref{thm:scr0msnegdepunanchored} to show that $(H_n(A)-\mu^2_A) \le 0$.
\end{proof}

%%%%%%%%%%%%%%%%%%%%%%%%%%%%%%%%%%%%%%%%%%%
\section{The scrambling advantage}
\label{sec:scrambling}

%[Would be nice to show two things: 1) that the first $n$ points of a scrambled $(0,s)$-sequence are also pairwise equidistributed over unachored boxes; 2) we can find an example of a digitally shifted $(0,m,s)$-net (or digitally shifted $P_n$ where $P_n$ are the first $n$ points of a $(0,s)$-sequence) that is far from being pairwise equidistributed, i.e., such that the variance of that digitally shifted point set is much worse than MC when integrating the indicator fct of that box. ]

We now give an example showing the advantage of scrambling over a digital shift, which is a simpler randomization. It uses a point set $P_n$ with $C_b(\BFk;P_n) \le 1$ %which, after randomization with a digital shift %, is 
such that
%, does not have a negative pairwise sampling index. In fact, we show that 
$P((\BFU,\BFV) \in A \times A) > {\rm Vol (A \times A)}$
for an anchored box $A$, for $(\BFU,\BFV)$ a pair of distinct points randomly chosen from the digitally shifted point set $\tilde{P}^{dig}_n$. So even the less restrictive condition ${\cal E}_{n,0}(\tilde{P}^{dig}_n)\le 0$ is not met. On the other hand, since $C_b(\BFk;P_n) \le 1$, Theorem 4.16 in \cite{WLD20} implies that %${\cal E}_n,{0}(\tilde{P}^{scr}_n) \le 0$ 
$P(\BFU \in A,\BFV \in A) \le {\rm Vol}^2(A)$ for $(\BFU,\BFV)$ randomly chosen from the scrambled point set $\tilde{P}_n$. 

\begin{example}
Consider the two-dimensional point set $P_n=\{(i/5,i/5),(i/5,((i+1) \bmod {5})/5),i=0,\ldots,4\}$. We first verify that $C_b(\BFk;P_n) \le 1$: this clearly holds for $\BFk = \mathbf{0}$. For $\BFk \in \{(1,0),(0,1)\}$, we have $C_b(\BFk;P_n) = 5 \times 1/9$. And for $\BFk$ with $|\BFk|\ge 2$ we have $C_b(\BFk;P_n) =0$. Now consider the box $A=[0,1/10)\times [0,2/5)$. Let us compute $P((\BFU,\BFV)\in A \times A)$, where $(\BFU,\BFV)$ is a pair of distinct points randomly chosen from $\tilde{P}^{dig}_n$, where $\tilde{P}_n^{dig} = P_n + \BFv$, where the addition is done digitwise and $\BFv \sim U(0,1)^2$. Let $\BFv = (v_1,v_2)$ with $v_j = 0.v_{j,1}v_{j,2}\ldots$, $j=1,2$. Then we see that among the $5^2$  possibilities for $(v_{1,1},v_{2,1})$, one point from $\tilde{P}^{dig}_n$ will be in the square $[0,1/5)\times[0,1/5)$ and one in the square $[0,1/5) \times [1/5,2/5)$ if and only if $v_{1,1}=v_{2,1}$, which happens with probability 1/5. Given that this happens, then it should also be clear that both points in that pair will be in $A$ if and only if $(0.0v_{1,2}v_{1,3}\ldots,0.0v_{2,2}v_{2,3}\ldots)\in [0,1/10) \times [0,1/5)$, which happens with probability 1/2 since $\BFv\sim U(0,1)^2$. Putting this all together, we get
$
P((\BFU,\BFV) \in A \times A) = \frac{1}{5}\frac{1}{45}\frac{1}{2}=\frac{1}{450},
$
where the fraction $1/45$ corresponds to the probability of choosing the pair of points falling in the squares (0,0) and (0,1) among the 45 different (unordered) pairs. Since ${\rm Vol}(A) = 1/25$, we have that $P((\BFU,\BFV) \in A \times A)
> ({\rm Vol}(A))^2 = 1/625$. 
\end{example}

\section{Future work}
\label{sec:conc}
In this paper, we have introduced a measure of uniformity for randomized QMC point sets that compares them to random sampling. This \emph{pairwise sampling dependence index} was shown to be no larger than 0 for scrambled $(0,m,s)-$nets, thus extending from anchored boxes to unanchored boxes the main result from \cite{WLD20}. For future work, we plan to try to extend our proof to the first $n$ points of a scrambled $(0,s)-$sequence. We also plan to explore how this result can lead to new bounds for the variance of scrambled $(0,m,s)-$nets in terms of the Monte Carlo variance for some functions.

%%%%
%\bibliography{Unanchored}

%\bibliographystyle{plain}

\begin{acknowledgement}
We thank the anonymous reviewers for their detailed comments. 
The first author thanks  NSERC for their support via grant \#238959. 
The second author wishes to acknowledge the support of the Austrian Science Fund (FWF):  Projects F5506-N26 and F5509-N26, which are parts of the Special Research Program ``Quasi-Monte Carlo Methods: Theory and Applications''.
\end{acknowledgement}

%%%%%%%%%%%%%%%%%%%%%%%%%%%%%%%%%%%%%%%%%%%%%%%%%%%%
\section*{Appendix: Proofs and Technical Lemmas}
\addcontentsline{toc}{section}{Appendix: Proofs and Technical Lemmas}

\label{sec:proofs}

We first prove results stated in the main part of the paper. These proofs make use of Lemmas \ref{lem:Ostr} to \ref{lem:weightedsums}, which are presented in the second part of the appendix.

\begin{proof}[of Lemma \ref{lem:vol}]
In what follows, we will use the notation $x_{\ell}$% 
%= \lfloor xb^{\ell} \rfloor$
\hjaspar{does this give you what you want? I do not think that this is the $\ell$th digit, it is the first $\ell$ digits as an integer.} to represent the $\ell$th digit in the base $b$ representation of $x \in [0,1)$, i.e., $x=\sum_{\ell\ge 1} x_{\ell} b^{-\ell}$, and the corresponding notation $x = 0.x_1 x_2 x_3 \ldots$. 
%We also use the notation $\tilde{V}_j = \sum_{i=j}^{\infty} V_i$ introduced in \cite{Lemieux2017}.

First, we decompose $A$ into three parts as $A_1 = [hb^{-r+1}+gb^{-r}-z,hb^{-r+1}+gb^{-r})$, $A_2 = [hb^{-r+1}+Gb^{-r}, hb^{-r+1}+Gb^{-r}+Z)$, $A_3 = [hb^{-r+1}+gb^{-r},hb^{-r+1}+Gb^{-r})$.
Hence we have
\hchristiane{Check notation for $V_i(A_1,A_2)$ wrt comment \# 17 of Ref 2}
\begin{equation}
\label{eq:decompVi}
V_i(A \times A) = \sum_{\ell=1}^3 V_i(A_{\ell} \times A_{\ell}) + 2 (V_i(A_1 \times A_2)+V_i(A_1\times A_3)+V_i(A_2\times A_3)) , \qquad i \ge 0.
\end{equation}
\hchristiane{We used $A_j$ to denote the right endpoint of the unanchored box, so maybe these $A_1,A_2,A_3$ should be denoted differently.}

Since $A_1$ and $A_2$ are both completely contained in the respective intervals $[hb^{-r+1}+(g-1)b^{-r},hb^{-r+1}+gb^{-r})$ and $[hb^{-r+1}+Gb^{-r},hb^{-r+1}+(G+1)b^{-r})$, any $x$ in $A_1$ is of the form $0.h_1\ldots h_{r-1}(g-1) x_{r+1}x_{r+2} \ldots$.
%, where the first $r-1$ digits are 0.
Similarly, $y \in A_2$ is of the form $0.h_1\ldots h_{r-1} (G) y_{r+1}y_{r+2} \ldots$.
\hjaspar{Similarly we need $0.\beta_1\ldots \beta_{r-1} (G) y_{r+1}y_{r+2} \ldots$} 
On the other hand, for $z \in A_3$ we have that $z_i=h_i$
\hjaspar{$z_i=\beta_i$} for $i \le r-1$, $z_r \in \{g,\ldots,G-1\}$, and $z_{\ell} \ge 0$ for $\ell>r$. From this we infer:
%\begin{enumerate}
%\item $V_i(A_j,A_j)=0$ for $i=0,\ldots,r-2, j=1,2,3$ (this implies $V_i(A,A)=0$ for $i \le r-2$, i.e., (i) holds);
%\item $V_{r-1}(A_j,A_j) = 0$ for $j=1,2$;
%\item $V_{r-1}(A_3,A_3) = \frac{(G-g)(G-g-1)}{b^{2r}}$;
%\item $V_{r-1}(A_i,A_j) = \Vol(A_i) \Vol(A_j)$ for $i \neq j$ because $\gamma_b(x,y)=r-1$ for any $x \in A_i$ and $y\in A_j$.
%\end{enumerate}
\hjaspar{I think the following is more clear than the list above.}
\hchristiane{Ok, but for second item I think it's $i=0,\dots,r-2$}
\hjaspar{I changed it}
\hchristiane{Ref 2 says doesn't see where we show $V_i(A_1,A_2) = 0$. Not sure we said we were proving this, i.e., only for $i=0,\ldots,r-2$. Could add that observation in item 3 below, since we're already noticing they have exactly $r-1$ initial digits in common}
\begin{enumerate}
    \item No pair of points from $A_1$ or $A_2$ can have less than $r$ initial common digits, thus $V_i(A_{\ell}\times A_{\ell})=0$ for $i=0,\dots,r-1$ and $\ell=1,2$.
    \item No pair of points from $A_3$ can have less than \hchristiane{$r-2$ or $r-1$?} {$r-1$} initial common digits, thus $V_i(A_3\times A_3)=0$ for $i=0,\dots,r-2$.
    \item A pair of points from any two of the following subsets: 
    $
    A_1,\, A_2,\, {[h b^{-r+1}+\beta b^{-r},}$ ${h b^{-r+1}+(\beta+1)b^{-r})}\subseteq A_3,
    $
    where ${\beta}=g,\dots, G-1$ has exactly $r-1$ initial common digits, thus $V_{r-1}(A_j\times A_{\ell})=\vol(A_j)\vol(A_{\ell})$ for $j\ne \ell$, $V_{r-1}(A_3 \times A_3)=(G-g)(G-g-1) b^{-2r}$
    {and $V_i(A_j \times A_{\ell}) = 0$ for $i \ge r, j \neq \ell$}.
\end{enumerate}
Note that this implies that $\tilde{V}_r(A_i \times A_i) = \Vol^2(A_i)$ for $i=1,2$, and (using item (3)) $\tilde{V}_{r}(A_3\times A_3) =\Vol^2(A_3) -V_{r-1}(A_3\times A_3)=  (G-g)b^{-2r}$.

%Property (i) in the statement of the Lemma is thus proved. 
The above statements also allow us to simplify \eqref{eq:decompVi}
as follows:
\begin{align}
V_{r-1}(A\times A) &= V_{r-1}(A_3\times A_3) + 2\sum_{1 \le i < \ell 
\le 3} V_{r-1}(A_i\times A_{\ell})\notag \\
&= V_{r-1}(A_3\times A_3) + 2 \sum_{1 \le i < \ell 
\le 3} \Vol(A_i)\Vol(A_{\ell}) \label{eq:Vrmoins1} \\
V_i(A\times A) &= \sum_{\ell=1}^3 V_i(A_{\ell}\times A_{\ell}) \quad i \ge r. \label{eq:Vipgr}
\end{align}

To prove (ii), consider the mappings $\varphi_j:[0,1) \rightarrow [0,1)$, $1\le j \le 3$ defined as:
\begin{align*}
\varphi_1(hb^{-r+1}+gb^{-r}-x) &= 1-x, \qquad 0 \le x < b^{-r}\\
\varphi_2(hb^{-r+1}+Gb^{-r}+x) &= x, \qquad 0 \le x < b^{-r}\\
%\end{align*}
%\begin{align*}
\varphi_3(hb^{-r+1}+gb^{-r}+x) &= x, \qquad 0 \le x <(G-g)b^{-r}.
\end{align*}\hjaspar{Same issue as before, we need $\beta b^{r-1}$}
All three are isometric mappings and such that
%\begin{align*}
$\varphi_1(A_1) = [1-z,1)$, 
$\varphi_2(A_2) = [0,Z)$, and 
$\varphi_3(A_3) = [0,(G-g)b^{-r})$.
%\end{align*}
Also, since $\varphi_j$ simply amounts to changing the first $r$ digits of a point in $A_j$ (and applies the same change to all points in $A_j$), it implies
\[
\gamma_b(\varphi_j(\nu_{j,\ell}),\varphi_j(\nu_{j,h})) = \gamma_b(\nu_{j,\ell},\nu_{j,h}), \qquad j=1,2,3,
\]
where $\nu_{1,\ell} =hb^{-r+1}+gb^{-r}-x,\nu_{1,h} = hb^{-r+1}+gb^{-r}-y$, $\nu_{2,\ell} =hb^{-r+1}+Gb^{-r}+x,\nu_{2,h} = hb^{-r+1}+Gb^{-r}+y$, and $\nu_{3,\ell} =hb^{-r+1}+gb^{-r}+w,\nu_{3,h} = hb^{-r+1}+gb^{-r}+z$.
%\begin{align*}
%\gamma_b(\varphi_1(hb^{-r+1}+gb^{-r}-x),\varphi_1(hb^{-r+1}+gb^{-r}-y))) &= \gamma_b(hb^{-r+1}+gb^{-r}-x,hb^{-r+1}+gb^{-r}-y) \\
%\gamma_b(\varphi_2(hb^{-r+1}+Gb^{-r}+x),\varphi_2(hb^{-r+1}+Gb^{-r}+y))) &= \gamma_b(hb^{-r+1}+Gb^{-r}+x,hb^{-r+1}+Gb^{-r}+y) \\
%\gamma_b(\varphi_3(hb^{-r+1}+gb^{-r}+w),\varphi_3(hb^{-r+1}+gb^{-r}+z))) &= \gamma_b(hb^{-r+1}+gb^{-r}+w,hb^{-r+1}+gb^{-r}+z). 
%\end{align*}
Therefore
\begin{align*}
V_i(A_1\times A_1)\hjaspar{Remove the subscript $i$?}
\hchristiane{I added $i$ on the RHS}
&= V_i(\varphi_1(A_1) \times \varphi_1(A_1)) = V_i([1-z,1) \times [1-z,1))  = V_i([0,z),[0,z))\\
V_i(A_2 \times A_2) &= V_i(\varphi_2(A_2) \times \varphi_2(A_2)) = V_i([0,Z) \times [0,Z)) \\
V_i(A_3 \times A_3) &= V_i(\varphi_3(A_3) \times\varphi_2(A_3)) = V_i([0,(G-g)b^{-r}) \times [0,(G-g)b^{-r})). 
\end{align*}
These intervals, being anchored at the origin,  satisfy the assumptions of Lemma 2.6 from \cite{WLD20}, which implies 
$bV_{i+1}(A_j \times A_j) -V_i(A_j \times A_j) \ge 0$ for $j=1,2,3$ and $i \ge 0$. Combining this with \eqref{eq:Vipgr}, property (ii) in the statement of Lemma \ref{lem:vol} is established.

For (iii), let us first assume $g=G$, and thus $A_3 = \emptyset$. Then, using \eqref{eq:Vrmoins1}, we get
$V_{r-1}(A \times A) = 2\Vol(A_1 \times A_2).
$
Furthermore, $\tilde{V}_r(A,A)   = \tilde{V}_r(A_1 \times A_1) + \tilde{V}_r(A_2 \times A_2) = \Vol^2(A_1) + \Vol^2(A_2)$. Since 
  $2\Vol(A_1 \times A_2) \le \Vol^2(A_1) + \Vol^2(A_2)$, (iii) is proved.
  
  Now assume $g<G$. In this case, we need to further refine $A_1$ and $A_2$ as: % and write them as
  \begin{align*}
  A_1 &= [hb^{-r+1}+gb^{-r}-db^{-(r+1)}-f, hb^{-r+1}+gb^{-r})\\
  A_2 &= [hb^{-r+1}+Gb^{-r},hb^{-r+1}+Gb^{-r}+Db^{-(r+1)}+F),
  \end{align*}
  where $0 \le d,D \le b-1$, $f,F \in [0,b^{-(r+1)})$. {Using \eqref{eq:Vrmoins1} and \eqref{eq:Vipgr}}, we then write
  \begin{align*}
  V_{r-1}(A \times A) &= V_{r-1}(A_3 \times A_3) + 2
  \sum_{1 \le i < \ell \le 3} \Vol(A_i)\Vol(A_{\ell}) \\
  &= \frac{(G-g)(G-g-1)}{b^{2r}}+2  \left( \left( \frac{d}{b^{r+1}}+f\right) \left( \frac{D}{b^{r+1}}+F\right)+ \frac{G-g}{b^r}\left(\frac{d}{b^{r+1}}+f\right) \right.\\
  &\left.+\frac{G-g}{b^r} \left( \frac{D}{b^{r+1}}+F\right)
    \right)
        \end{align*}
    \begin{align*}
  V_r(A \times A) &= \sum_{\ell=1}^3 V_r(A_{\ell} \times A_{\ell}) 
 = V_r(A_1 \times A_1)+V_r(A_2 \times A_2) + \frac{G-g}{b^{2r}} \frac{b-1}{b} \\
  &= \frac{(d-1)d}{b^{2(r+1)}} + \frac{2fd}{b^{r+1}}+\frac{(D-1)D}{b^{2(r+1)}} + \frac{2FD}{b^{r+1}} + \frac{G-g}{b^{2r}} \frac{b-1}{b} \\
  \tilde{V}_r(A \times A) &= \left(\frac{d}{b^{r+1}}+f\right)^2 +  \left(\frac{D}{b^{r+1}}+F\right)^2 + \frac{G-g}{b^{2r}}.
  \end{align*}
  The last equality for $V_r(A \times A)$ is obtained by observing that for $(x,y)$ to be in $V_r(A_1 \times A_1)$, either (i) $x=0.h_1\ldots h_{r-1}(g-1)d_1\ldots$ and $y=0.h_1\ldots h_{r-1}(g-1)d_2\ldots$ with $d_1 \neq d_2 \in \{0,\ldots,d-1\}$, or (ii) one of them is of the form $z_1 + (0.h_1\ldots h_{r-1}(g-1)d) $
  with $z_1 \in [0,f)$ and the other is of the form $0.h_1\ldots h_{r-1}(g-1)d_1\ldots$ with $d_1 \in \{0,\ldots,d-1\}$. Case (i) contributes a volume of size $(d-1)d b^{-2(r+1)}$ and case (ii) contributes $2fdb^{-(r+1)}$.
  A similar argument can be used to derive $V_r(A_2 \times A_2)$.
\hjaspar{yes please}
  \hchristiane{We may consider adding a picture to explain this further}
  
  Therefore $V_{r-1}(A \times A) -\frac{b(b-2)}{b-1}V_r(A \times A) \le \tilde{V}_r(A \times A)$ holds if 
  \begin{align}
 &\frac{(G-g)(G-g-1)}{b^{2r}}+2 \left(\frac{G-g}{b^r} \left( \frac{d+D}{b^{r+1}}+f+F\right)\right)  +  2\left( \frac{d}{b^{r+1}}+f\right) \left( \frac{D}{b^{r+1}}+F\right) \notag \\
  -&\frac{b(b-2)}{b-1}
  \left( \frac{(d-1)d}{b^{2(r+1)}} + \frac{2fd}{b^{r+1}}+\frac{(D-1)D}{b^{2(r+1)}} + \frac{2FD}{b^{r+1}}+ \frac{G-g}{b^{2r}} \frac{b-1}{b}\right)  \notag \\
  \le&  \left(\frac{d}{b^{r+1}}+f\right)^2 +  \left(\frac{D}{b^{r+1}}+F\right)^2 + \frac{G-g}{b^{2r}}.
  \label{ineq:gGdDfF}
  \end{align}
  Since 
  \[
   2\left( \frac{d}{b^{r+1}}+f\right) \left( \frac{D}{b^{r+1}}+F\right) 
   \le  \left(\frac{d}{b^{r+1}}+f\right)^2 +  \left(\frac{D}{b^{r+1}}+F\right)^2
  \]
  it means that to prove \eqref{ineq:gGdDfF} it is sufficient to show that
  \begin{align*}
 & \frac{(G-g)(G-g-1)}{b^{2r}} - (b-2) \frac{G-g}{b^{2r}}+2 \left(\frac{G-g}{b^r} \left( \frac{d+D}{b^{r+1}}+f+F\right)\right) \\
&   -\frac{b(b-2)}{b-1}
  \left( \frac{(d-1)d}{b^{2(r+1)}} + \frac{2fd}{b^{r+1}}+\frac{(D-1)D}{b^{2(r+1)}} + \frac{2FD}{b^{r+1}} \right) 
  \le  \frac{G-g}{b^{2r}},
  \end{align*}
  or equivalently, that
    \begin{align}
 & -b^2 (G-g)(b-(G-g)) +2b (G-g) \left( d+D+b^{r+1}(f+F)\right) \notag \\
&   -\frac{b(b-2)}{b-1}
  \left( (d-1)d+ 2fdb^{r+1}+(D-1)D+ 2FDb^{r+1}\right) 
  \le  0. \label{eq:MainIneqiiiLemma02}
  \end{align}
  
Note that $G-g \le b-2$ by assumption.  We proceed by considering three cases:
  
  \noindent {\bf Case 1:} $G-g \le b-4$.
This implies $b-(G-g) \ge 4$ {(and thus $b \ge 4$)}. Also, to handle this case we use the fact that $0 \le fb^{r+1},Fb^{r+1} < 1$. 
By making appropriate substitutions for $f$ and $F$, we see that to prove \eqref{eq:MainIneqiiiLemma02} holds it is sufficient to show that
  \[
  -4b^2(G-g)+2b(G-g)(d+D+2) - \frac{b(b-2)}{b-1}(d(d-1)+D(D-1)) \le 0
  \]
  which holds because $d+D+2 \le 2b$.
  
\noindent {\bf Case 2: $G-g = b-3$}
{First note that this implies $b \ge 3$. Next, }
we replace $G-g$ with $b-3$ in  \eqref{eq:MainIneqiiiLemma02} and divide each term by $b$. 
For this case, we can use the bound $0 \le fb^{r+1},Fb^{r+1} <1$ and by substituting appropriately, it means it is sufficient to show
%
% so that the inequality to prove is
%\begin{align*}
%-3b(b-3)+2(b-3)(d+D+(f+F)b^{r+1})-\frac{b-2}{b-1}(d(d-1)+D(D-1)+2b^{r+1}(fd+FD)) \le 0.
%\end{align*}
\[
-3b(b-3)+2(b-3)(d+D+2)-\frac{b-2}{b-1}(d(d-1)+D(D-1)) \le 0.
\]
We view the LHS as the sum of two quadratic polynomials, $p(d)$ and $p(D)$, and thus argue it is sufficient to show that
\[
p(d) := \frac{-(b-2)}{b-1} d^2 + d \left(2(b-3) +\frac{b-2}{b-1}  \right)
-(b-3)(3b/2-2) \le 0.
\]
\hjaspar{What happens to the second polynomial? Why can we ignore it when proving the inequality? By my calculations the inequality is not $p(d)+p(D)\leq 0$. (same problem as the case below)}
\hchristiane{There was a mistake in the constant of $p(d)$, which should be $-(b-3)(3b/2-2)$ i.e., need to divide by 2 so it's the same in $p(d)$ and $p(D)$. That is, we rewrite the expression in $d$ and $D$ as being $p(d)+p(D)$ where $p()$ is a polynomial of degree 2. So if we can prove $p(d) \le 0$ for all $d$, then it means $p(d)+p(D) \le 0$. It's fixed now.}
We will show this holds by finding the value $d_{max}$ of $d$ that maximizes $p(d)$ and show that $p(d_{max}) \le 0$.
We have that
\[
p'(d) = -2d\frac{b-2}{b-1} + 2(b-3) +\frac{b-2}{b-1}.
\]
Therefore
\[
d_{max} = \left(2(b-3) +\frac{b-2}{b-1}  \right) \frac{b-1}{2(b-2)} = \frac{(b-3)(b-1)}{b-2} + \frac{1}{2}.
\]
Hence $d_{max} \in (b-2.5,b-1.5)$. Thus it is sufficient to show $p(b-2) \le 0$.
Now,
\[
p(b-2) = -\frac{(b-2)^3}{b-1} + (b-2)\left(2(b-3)+\frac{b-2}{b-1}\right) -(b-3)(3b/2-2) 
\]
therefore
\begin{align*}
(b-1)p(b-2)&=-(b-2)^3+2(b-1)(b-2)(b-3)+(b-2)^2-\left(\frac{3b}{2}-2\right)(b-3)(b-1) \\
%&= (b/2-2)(b-3)(b-1)+(b-2)^2(3-b) \\
&=(3-b)(b^2-3b+4)/2=(3-b)(b(b-3)+4)/2\le 0 
\end{align*}\hjaspar{Is there some trick to the last equality? Or do you have to expand both sides and verify the equality?}
\hchristiane{This has changed so should be more obvious now.}
since $b \ge 3$.

\noindent {\bf Case 3: $G-g = b-2$}

In this case \eqref{eq:MainIneqiiiLemma02} becomes
\begin{align*}
&-2b^2(b-2)+2b(b-2)(d+D+(f+F)b^{r+1})\\ 
&- \frac{b(b-2)}{b-1} \left(d(d-1)+D(D-1) + 2b^{r+1} (fd+FD)\right)  \le 0 \\
\Leftrightarrow & -2b+2(d+D+(f+F)b^{r+1})  - \frac{1}{b-1} \left(d(d-1)+D(D-1) + 2b^{r+1} (fd+FD)\right)  \le 0.
\end{align*}
As in the case $G-g=b-3$, we argue it is sufficient to show each of the quadratic polynomials in $d$ and $D$ on the LHS (which are the same) is bounded from above by 0. That is, we need to show
\[
q(d) := \frac{-d^2}{b-1} + d\left(2+\frac{1}{b-1}-\frac{2fb^{r+1}}{b-1} \right) -(b-2fb^{r+1}) \le 0.
\]
\hjaspar{For the polynomials to be the same the last term needs to be replaced with $-b+2fb^{r+1}$.}
\hchristiane{Fixed. Means we get exactly $q(b-1)=0$ below.} 
Now
\[
q'(d) = \frac{-2d}{b-1} +2 + \frac{1}{b-1}-\frac{2fb^{r+1}}{b-1}
\]
and thus $d_{max} = \left(2 + \frac{1}{b-1}-\frac{2fb^{r+1}}{b-1} \right) \frac{b-1}{2} = b-0.5-fb^{r+1}$, which implies
$d_{max} \in (b-1.5, b-0.5)$. Thus it is sufficient to show $q(b-1) \le 0$. We have that
\begin{align*}
q(b-1) &= -(b-1)+(b-1)\left(2 + \frac{1}{b-1}-\frac{2fb^{r+1}}{b-1}\right)-(b-2fb^{r+1})  \\
&=(b-1)+1-2fb^{r+1}-b+2fb^{r+1} =  0
\end{align*}
as required.
\end{proof}

\begin{proof}[of Theorem \ref{thm:decompA}]
To simplify the notation, we define 
$Z_k := b^{2k}V(1_k \times I_k)$ and ({$
W_k^{(r-1)} := (b^{2(k+r+1)}/4)V(Y_k^{(r-1)} \times Y_k^{(r-1)})$} to be  the (normalized) volume vectors of $k-$elementary intervals and elementary unanchored $(k,r-1)-$intervals, respectively. {Since the coordinates of each of these vectors are positive and sum to one, $V(A\times A)=\sum_{k\geq 0}(\alpha_k+\tau_k)$ follows immediately from the last equality in the statement of the theorem.} 

Based on the definition of $Y_k^{({r-1})}$, the vectors $W_k^{({r-1})}$
satisfy, for $i\ge 0$, {$r \ge 1$},
\hchristiane{Alternatively to what I added in blue, could switch 2nd and 3rd case and then say it's 0 otherwise}
\[
W_{i,k}^{({r-1})} = \begin{cases}
1/2 & \mbox{if }i=r-1 \\
(b-1)/2b^{i-(k+{r})} & \mbox{if } i \ge k+{r+1}\\
0 & \mbox{{otherwise},}%\mbox{if } {r-1<i<k+r+1} \mbox{  or } i<r-1,
\end{cases}
\]
while
\[
Z_{i,k} = \begin{cases}
0 & \mbox{if } i<k \\
(b-1)/b^{i-k+1} & \mbox{if } i \ge k.
\end{cases}
\]
Note that 
{
\begin{align}
    W_{r-1,k}^{(r-1)} &= 1/2 \mbox{ for all }k \ge 0 \label{Wdkdonehalf}\\
    W_{i,k}^{(r-1)} &= Z_{i,k+r+1}/2 \mbox{ for }i \ge k+r+1.\label{eq:WikdvsZik}
\end{align}
}

There are two cases to consider. First, if $V_{r-1}(A \times A) \le bV_r(A \times A)$, then based on Properties (i) and (ii) from Lemma \ref{lem:vol}  we can decompose $V(A \times A)$ solely with the $Z_k$'s \hjaspar{Better to use $Z_k$.}, i.e., we set $\alpha_k = 0$ for all $k$ and 
\[
\tau_k = \frac{bV_k(A \times A)-V_{k-1}(A \times A)}{b-1} \quad k \ge r -1
\]
and $\tau_k = 0$ if $0\le k\le r-2$. {Note that $A=[0,1)$ fits into this first case.}

Second, if $V_{r-1}(A \times A) \ge bV_r(A \times A)$, then we first decompose the vector $\sum_{k=r}^{\infty}V_k(A \times A) \BFe_k$, where $\BFe_k$ is  a (canonical) vector of zeros with a 1 in position $k$, (note that this agrees with the vector $V(A \times A)$ everywhere except on index $r-1$, where it has a 0 instead of $V_{r-1}(A \times A)$) as
\[
\sum_{k=r}^{\infty}V_k(A \times A) \BFe_k = \sum_{k=r}^{\infty} \bar{\tau}_k Z_k
\]
with $\bar{\tau}_r = bV_r(A \times A)/(b-1) \ge 0$, $\bar{\tau}_k=0$ if $k <r$ (from Part (i) of Lemma \ref{lem:vol}), and 
\[
\bar{\tau}_k =  \frac{bV_k(A \times A)-V_{k-1}(A \times A)}{b-1} \quad k \ge r +1.
\]
From Lemma \ref{lem:vol} we know $\bar{\tau}_k \ge 0$ for $k \ge r+1$.
Note that $b\bar{\tau}_r Z_{r-1}- \bar{\tau}_r Z_r = bV_r(A \times A) \BFe_{r-1}$, i.e., $b\bar{\tau}_r Z_{r-1}$ agrees with $\bar{\tau}_r Z_r$ everywhere except on index $r-1$. Therefore
\[
V(A \times A) - b\bar{\tau}_r Z_{r-1}  - \sum_{k=r+1}^{\infty} \bar{\tau}_k Z_k = (V_{r-1}(A \times A)-bV_r(A \times A)) \BFe_{r-1}.
\]
Hence the $\sum_{k = 0}^{\infty} \alpha_k W_k^{(r-1)} $ part of the decomposition is only needed to decompose $V_{r-1}(A \times A)-bV_r(A \times A)$.
We claim there exists $\bar{\alpha}_k \ge 0$, $k \ge r+1$ such that
\begin{equation}
\label{eq:alphaksurdeux}
V_{r-1}(A \times A)-bV_r(A \times A) = \sum_{k=r+1}^{\infty} \bar{\alpha}_k/2,
\end{equation}
and such that $\bar{\alpha}_k/2 \le \bar{\tau}_k$ for $k \ge r+1$. This can be seen using Part (iii) of Lemma \ref{lem:vol}. Indeed, to ensure the existence of these $\bar{\alpha}_k$'s, we need to prove that
\begin{equation}
\label{eq:Vr1bvrlesumtau}
V_{r-1}(A \times A)-bV_r(A \times A)  \le \sum_{k=r+1} \bar{\tau}_k.
\end{equation}
Now,
\begin{align*}
 \sum_{k=r+1} \bar{\tau}_k &= \sum_{k=r+1}^{\infty} \frac{bV_k(A \times A)-V_{k-1}(A \times A)}{b-1} \\
 &= \tilde{V}_{r+1}(A \times A) - \frac{V_r(A \times A)}{b-1} = \tilde{V}_r(A \times A) -V_r(A \times A) -\frac{V_r(A \times A)}{b-1}.
\end{align*}
Therefore, \eqref{eq:Vr1bvrlesumtau} holds if and only if
\[
V_{r-1}(A \times A)-bV_r(A \times A) \le  \tilde{V}_r(A \times A) -V_r(A \times A) -\frac{V_r(A \times A)}{b-1}  %\Leftrightarrow 
\]
which holds if and only if
$
V_{r-1}(A \times A) \le 
 \tilde{V}_r(A \times A) +  \frac{b(b-2)}{b-1}V_r(A \times A),
$
which is precisely what Part (iii) of Lemma \ref{lem:vol} shows. Having proved the existence of non-negative coefficients $\bar{\alpha}_k$ satisfying \eqref{eq:alphaksurdeux} implies we can write
$
V_{r-1}(A \times A)-bV_r(A \times A) = \sum_{k=r+1}^{\infty} \bar{\alpha}_k W_{r-1,k-(r+1)}^{(r-1)}
$
by using \eqref{Wdkdonehalf}. Hence all that is left to do is to find the combination of $Z_k$'s that can cancel out $\sum_{k=r+1}^{\infty} \bar{\alpha}_k W_{i,k-(r+1)}^{(r-1)}$ for $i \ge r+1$ (we can ignore the case $i=r$ because $W_{r,k-(r+1)}^{(r-1)}=0$ for all $k \ge r+1$). This is done by using \eqref{eq:WikdvsZik}, which implies that
\[
\sum_{k=r+1}^{\infty} \bar{\alpha}_k W_{i,k-(r+1)}^{(r-1)}
=\sum_{k=r+1}^{\infty} (\bar{\alpha}_k/2) Z_{i,k}.
\]
Hence the final decomposition is given by
\begin{align*}
V(A \times A) &= \sum_{k=r+1}^{\infty} \bar{\alpha}_k W_{k-(r+1)}^{(r-1)} +  b\bar{\tau}_r Z_{r-1}  + \sum_{k=r+1}^{\infty} (\bar{\tau}_k-\bar{\alpha}_k/2) Z_k\\
&= \sum_{k=0}^{\infty} \alpha_k W_k^{(r-1)} + \sum_{k=0}^{\infty} \tau_k Z_k
\end{align*}
with $\alpha_k = \bar{\alpha}_{k+(r+1)}$, $k \ge 0$, 
$\tau_k = \bar{\tau}_k - \bar{\alpha}_k/2$, $k \ge r+1$,
$\tau_{r-1} = b \bar{\tau}_r$ and $\tau_k =0$ for $0 \le k \le r-2, k=r$.
\end{proof}

\begin{proof}[of Lemma \ref{lem:volcondprob}]
(1) From the definition of $D(\BFk,\BFd,J)$, we have that
\[
{\rm Vol}(D(\BFk,\BFd,J)) = b^{-2|\BFk|_{J^c}}2^{2|J|}b^{-2(|\BFk+\BFd+2|_J)}=
2^{2|J|} b^{-2(|\BFk|+|\BFd+2|_J)}.
\]
(2) Similarly, from the definition of $F_1(\BFk,\BFd,J,I)$ we get 
\[
{\rm Vol}(F_1(\BFk,\BFd,J,I)) = b^{-|\BFk|_{J^c}}b^{-(|\BFk+\BFd+2|_J)} = b^{-(|\BFk|+|\BFd+2|_J)}.
\]
(3) This conditional probability is given by $\eta/n-1$, where $\eta$ is the number of points $\BFu_{\ell}$ with $\ell \neq i$ that are in $F_2$ if $\BFu_i \in F_1$.
Hence for $j \in I$ we must have $\gamma_b(u_{i,j},u_{\ell,j}) \ge k_j+d_j+2$; for $j \in I^c$ we must have $\gamma_b(u_{i,j},u_{\ell,j}) \ge k_j$.
For $j \in J \cap I^c$, the requirement that $u_{i,j} \in F_1$ means $u_{\ell,j}$ must satisfy: %the following conditions:
\begin{itemize}
    \item[(a)] it must have the same first $d_j$ digits as $u_{i,j}$;
    \item[(b)] its $(d_j+1)$th digit must be 1 (while the $(d_j+1)$th digit of $u_{i,j}$ is 0);
    \item[(c)] the digits $u_{\ell,j,r}$ for $d_j+2 \le r d_j+k_j+2$ must be 0
\end{itemize}
If we only had to satisfy requirement (a), then we would have $\eta = m_b(\BFk,\BFd,2,J,I)$. However, the requirements (b) and (c) imply 
\[
\eta =  m_b(\BFk,\BFd,2,J,I) 
\prod_{j \in J \cap I^c} \frac{1}{b-1}
\frac{1}{b^{k_j+1}},
\]
where the term $1/(b-1)$ handles restriction (b)
while the term $b^{-k_j-1}$ handles restriction (c).
Therefore
\begin{align*}
P(\BFV \in F_2(\BFk,\BFd,J,I)|\BFU \in F_1(\BFk,\BFd,J,I))
 &= \frac{m_b(\BFk,\BFd,2,J,I)}{n-1} \left(\frac{1}{b-1}\right)^{|J|-|I|}
 \frac{1}{b^{|\BFk+1|_{J \cap I^c}}}\\
&= \frac{m_b(\BFk,\BFd,2,J,I)}{n-1} \frac{(b-1)^{|I|-|J|}}{b^{|\BFk|_{J \cap I^c}+|J|-|I|}}.
\end{align*}
(4) The decomposition  $\cup_{K,I \subseteq J} F(\BFk,\BFd,J,I,K)$ is obtained by expanding each $Y_{k_j}^{(d_j)}$ as $Y_{k_j,1}^{(d_j)} \cup Y_{k_j,2}^{(d_j)}$. Then, we need to prove that 
for $K_1,K_2 \subseteq J$, 
\begin{equation}
\label{eq:puvk1k2}
    P(\BFU,\BFV) \in F(\BFk,\BFd,J,I,K_1) = P(\BFU,\BFV) \in F(\BFk,\BFd,J,I,K_2). 
\end{equation}
Noting that the equality \eqref{eq:HAinsidelemma} can be generalized to $P((\BFU,\BFV) \in {\cal R})
=\sum_{\BFi \ge \mathbf{0}} \psi_{\BFi}V_{\BFi}({\cal R})$, it is clear that to prove \eqref{eq:puvk1k2}, it is sufficient to show that 
the volume vectors corresponding to $F(\BFk,\BFd,J,I,K_1)$ and $F(\BFk,\BFd,J,I,K_2)$ are equal. To do so, since each entry $V_{\BFi}({\cal R}) = \prod_{j=1}^s V_{i_j}({\cal R}_j)$ (where for ${\cal R} = {\cal R}_1 \times {\cal R}_2$ we write ${\cal R}_j = {\cal R}_{1,j} \times {\cal R}_{2,j}$),  it is sufficient to show that for fixed $k$ and $d$, $V_{11}:= V(Y_{k,1}^{(d)} \times Y_{k,1}^{(d)}) =  V(Y_{k,2}^{(d)} \times Y_{k,2}^{(d)})=: V_{2,2}$
and $V_{12} := V(Y_{k,1}^{(d)} \times Y_{k,2}^{(d)})=
V(Y_{k,2}^{(d)} \times Y_{k,1}^{(d)})=: V_{2,1}$.
But this follows from an argument similar to the one used in the proof of Lemma 4.14 in \cite{WLD20}, which we adapt for our setup. First we introduce  the set
${\cal F} = \{ab^{-k}: a\in \mathbb{Z},k \in \mathbb{N}\} \subseteq \mathbb{R}$ which has Lebesgue measure 0. Then, we argue that for $x,y \in (0,b^{-(k+d+2)}) \cap {\cal F}^c$, we have
\begin{align*}
    \gamma_b\left(\frac{1}{b^{d+1}}-x,\frac{1}{b^{d+1}}-y\right) = \gamma_b\left(\frac{1}{b^{d+1}}+x,\frac{1}{b^{d+1}}+y\right) \\
    \gamma_b\left(\frac{1}{b^{d+1}}-x,\frac{1}{b^{d+1}}+y\right) = \gamma_b\left(\frac{1}{b^{d+1}}+x,\frac{1}{b^{d+1}}-y\right).
\end{align*}
Therefore, for $(x,y) \in Y_{k}^{(d)} \times Y_k^{(d)} $, $D_i$
is, up to a set of measure 0, invariant under the transformation $(x,y) \mapsto \left(\frac{2}{b^{d+1}}-x,\frac{2}{b^{d+1}}-y\right)$. This transformation maps $Y_{k,1}^{(d)} \times Y_{k,1}^{(d)}$ to $Y_{k,2}^{(d)}  \times Y_{k,2}^{(d)}$ (and vice-versa) and $Y_{k,1}^{(d)} \times Y_{k,2}^{(d)}$ to $Y_{k,2}^{(d)}  \times Y_{k,1}^{(d)}$ (and vice-versa). Therefore
$V_{11}=V_{22}$, and
$V_{12} = V_{21}$, as required.
\end{proof}

\begin{proof}[of Lemma \ref{gjiloosebound}]
First, if $\ell = 2i+1$, then
\[
g_{i,j}(\ell) = \left(\frac{b}{b-1}\right)^{j-i-1}
\le \left(\frac{b}{b-1}\right)^{i+j-\ell}
\]
and
\[
\left(\frac{b}{b-1}\right)^{i+j-\ell}
= \left(\frac{b}{b-1}\right)^{i+1-j}
\]
so in this case we actually have an equality, i.e.,
\[
g_{i,j}(\ell)  =  \left(\frac{b}{b-1}\right)^{i+j-\ell}.
\]
If $2i< \ell < i+j$ and $\ell$ is even, then
\[
\left(\frac{b}{b-1}\right)^{i+j-\ell} \ge 1 = g_{j,i}(\ell)
\]
since $\ell < i+j$.
If $2i< \ell < i+j$ and $\ell$ is odd, then
we must show that 
%What is left to do is show that when $k>1$ is odd, then 
\begin{equation}
\label{eq:ineg1}
1 +  \frac{b^{j+i-\ell}}{(b-1)^{j-i}} \binom{j-i-1}{\ell-2i} \le  \left(\frac{b}{b-1}\right)^{i+j-\ell}.
\end{equation}
Now let $k = \ell-2i$ and $r=j-i$. This means $k$ is odd and $k>1$, and also $k < r \le j \le s$ which means $r \ge 4$.
Using this notation, \eqref{eq:ineg1} is equivalent to
%\[
%\frac{(b-1)^r}{b^{r-k}}+\binom{r-1}{k}
%\le (b-1)^k.
%\]
%Hence the inequality that needs to be established is
 \begin{equation}
 \label{eqn:bm1keven}
 b^k \left(\frac{b-1}{b}\right)^r +\binom{r-1}{k} \le (b-1)^k %\Leftrightarrow
 \mbox{ {and thus to} }
  \left(\frac{b-1}{b}\right)^{r-k} + \frac{\binom{r-1}{k}}{(b-1)^k} \le 1.
 \end{equation}
 %holds. 
 Now,
 \[
  \left(\frac{b-1}{b}\right)^{r-k} = \sum_{j=0}^{r-k}  \binom{r}{j} \left(\frac{-1}{b}\right)^{j}
 \]
and $\binom{r}{j}/b^j$ is decreasing with $j$, therefore
\[
\left(\frac{b-1}{b}\right)^{r-k} \le 1 - \frac{r-k}{b}+\frac{(r-k)(r-k-1)}{2b^2},
\]
because the condition that $k>1$ and $k$ is odd implies $k \ge 3$.
Therefore a sufficient condition for {the second inequality in} \eqref{eqn:bm1keven} to hold is if we have
\begin{align}
& 1-\frac{r-k}{b}+\frac{(r-k)(r-k-1)}{2b^2} +  \frac{(r-1)(r-2)\ldots (r-k)}{(b-1)(b-1)\ldots (b-1)} \frac{1}{k!} \le 1  \notag\\
%\Leftrightarrow 
\mbox{{which holds iff }}
& \frac{(r-k)(r-k-1)}{2b^2} +  \frac{(r-1)(r-2)\ldots (r-k)}{(b-1)(b-1)\ldots (b-1)} \frac{1}{k!} \le \frac{r-k}{b} \notag \\
%\Leftrightarrow 
\mbox{{which holds iff }}%
& \frac{r-k-1}{2b} + \frac{(r-1)(r-2)\ldots (r-k+1)}{(b-1)(b-1)\ldots (b-1)}  \frac{b}{b-1}\frac{1}{k!} \le 1.
\label{eq:triangle}
\end{align}
Now, 
\[
\frac{r-k-1}{2b} < \frac{1}{2} \Leftrightarrow r-k-1 < b,
\]
and the latter inequality holds since $b \ge s \ge r$ and $k \ge 3$. Also,
\[
\frac{r-k+1}{2(b-1)} \le  \frac{1}{2} \Leftrightarrow r-k+1 \le  b -1 \Leftrightarrow r+2-k \le b
\]
and the latter inequality holds since $k \ge 3$ and $b \ge s \ge r$. Therefore, for $k \ge 3$ the LHS of \eqref{eq:triangle}  is bounded (strictly) from above by
\[
\frac{1}{2} + \frac{1}{2}  \frac{(r-1)(r-2)\ldots (r-k+2)}{(b-1)(b-1)\ldots (b-1)}  \frac{b}{b-1} \frac{1}{k(k-1)\ldots 3} < 1
\]
because 
\[
\frac{(r-1)(r-2)\ldots (r-k+2)}{(b-1)(b-1)\ldots (b-1)} \le 1 
\]
since $b \ge s \ge r$ and 
\[
 \frac{b}{b-1} \frac{1}{k(k-1)\ldots 3} \le \frac{4}{3} \times \frac{1}{3} < 1
\]
since $b/(b-1)$ decreases with $b$, and the condition $1<k<r\le s$ with $k \ge 3$ means we can assume $b \ge s \ge r \ge 4$. This proves that  \eqref{eqn:bm1keven} holds.
\end{proof}

\begin{proof}[of Lemma \ref{lem:boundPsi}]
(i) Using Lemma \ref{lem:Mtilde} with $c=2$ (which implies $c|I|+|I^{*c}| = |J|+|I|$), we first consider the case where  $m-|\BFk|_{I^*}-|\BFd|_J \ge |I|+|J|$. In this case
\[
%\tilde{n}_{I^c}(\BFk+\mathbf{2}+\BFd,\BFd)  = 
m(\BFk,\BFd,2,J,I)  %= 
  %(b-1)^{|J \cap I^c|} b^{m- \bar{\nu}_I}
  = (b-1)^{|J|-|I|} b^{m-|\BFk|_{I^*}-|\BFd|_J-|J|-|I| }.
\]
Hence if $m-|\BFk|_{I^*}-|\BFd|_J \ge |I|+|J|$, then
%we can bound $\psi_m(\BFk,\BFd,J,I) $ by
\begin{align*}
%g(m,I,\BFk,\BFd) &=  
\psi_m(\BFk,\BFd,J,I) =\frac{(b-1)^{|J|-|I|}}{n-1} b^{m-|\BFk|_{I^*}- |\BFd|_J -|J|-|I| } b^{|\BFk|_{I^*}+|\BFd|_J +|J|+|I|}(b-1)^{|I|-|J|} 
= \frac{b^{m}}{b^m-1}.
\end{align*}
(ii) Next, again based on  Lemma \ref{lem:Mtilde}, we consider the case  %$|\BFk|_I+|\BFd|_{I}+2|I| +  %|\BFk|_{J^c}+|\BFd|_{J \cap I^c} +1<m < %|\BFk|_I+|\BFd|_{I}+2|I| +|\BFk|_{J^c}+|\BFd|_{J %\cap I^c}+|J|-|I|$ or equivalently, 
$2|I| +1 <m-|\BFk|_{I^*}-|\BFd|_J < |I|+|J|$.
First, if $m-|\BFk|_{I^c}-|\BFd|_J$ is odd and $m - |\BFk|_{I^c}-|\BFd|_J>2|I|+1$, then

%[I need to fix this and start here to separate into the two cases i) and ii) used below]
\begin{align}
%\tilde{n}_{I^c}(\BFk+\mathbf{2}+\mathbf{d},\BFd)
m(\BFk,\BFd,2,J,I)
%&\le 
%\left(b^{m-|\BFk|_{I^*}-|\BFd|_{I}-2|I|-|\BFd|_{J \cap I^c} } \left(\frac{b-1}{b}\right)^{|J|-|I|}
%+\binom{|J|-|I|-1}{m-|\BFk|_I-|\BFd|_I-2|I|-|\BF%k|_{J^c}-|\BFd|_{J \cap I^c} } \right)\notag\\
&\le  \left(b^{m-|\BFk|_{I^*}-|\BFd|_J-2|I| } \left(\frac{b-1}{b}\right)^{|J|-|I|}
+\binom{|J|-|I|-1}{m-|\BFk|_{I^*}-|\BFd|_J-2|I|}\right). \label{Mbd1}
%&\le (b-1)^{m-|\BFk|_{I^*}- |\BFd|_J-2|I| }. %\label{Mbd2}
\end{align}
%where 
%\[
%{\bf 1}_{m,I,\BFk,\BFd} = \begin{cases}
%1 & \mbox{ if } m-|\BFk|_I-|\BFd| \mbox{ is odd} \\
%0 & \mbox{otherwise}.
%\end{cases}
%\]
%April 10

Hence in that case %if $m - |\BFk|_I-|\BFd|-2|I|>1$ and is odd, then
\begin{align*}
 \psi_m(\BFk,\BFd,J,I)  &\le   \frac{1}{b^m-1} b^{|\BFk|_{I^*}+|\BFd|_J+|J|+|I|}(b-1)^{|I|-|J|}  
 \left(b^{m-|\BFk|_{I^*}-|\BFd|_J-2|I| } \left(\frac{b-1}{b}\right)^{|J|-|I|}\right.\\
&\left.+\binom{|J|-|I|-1}{m-|\BFk|_{I^*}-|\BFd|_J-2|I|} \right)\\
&= \frac{1}{b^m-1} \left(b^m +b^{|\BFk|_{I^*}+|\BFd|_J+|J|+| I |}(b-1)^{|I|-|J|}
\binom{|J|-|I|-1}{m-|\BFk|_{I^*}-|\BFd|_J-2|{I}|} \right) \\
&= \frac{b^m}{b^m-1}  \left( 1+\frac{b^{|J|+ |I| -(m -|\BFk|_{I^*}-|\BFd|_J)}}{(b-1)^{|J|-|I|} }
\binom{|J|-|I|-1}{m-|\BFk|_{I^*}-|\BFd|_J-2|{I}|} \right). 
%&=\frac{b^m}{b^m-1}  g_1(m,I,\BFk,\BFd)
\end{align*}
A similar calculation shows that  if $m - |\BFk|_{I^*}-|\BFd|_J$ is even then 
 \[
  \psi_m(\BFk,\BFd,J,I) \le  \frac{b^m}{b^m-1}.
\] 
If $m - |\BFk|_{I^c}-|\BFd|_J-2|I|=1$, then from  Lemma \ref{lem:Mtilde} we know that 
$m(\BFk,\BFd,2,J,I) =(b-1)$, and thus
 \begin{align*}
 \psi_m(\BFk,\BFd,J,I)  =  \frac{1}{b^m-1} b^{|\BFk|_{I^*}+|\BFd|_J+|J|+|I|}(b-1)^{|I|-|J|} (b-1) \\
 = \frac{b^m}{b^m-1} b^{|J|-|I|-1}(b-1)^{1+|I|-|J|}  = \frac{b^m}{b^m-1} \left(\frac{b}{b-1}\right)^{|J|-|I|-1}.
 \end{align*}
 Combining these three cases, we get that for $2|I| < m-|\BFk|_{I^*}-|\BFd|_J < |J|+|I|$,
 \[
  \psi_m(\BFk,\BFd,J,I)  =  \frac{b^m}{b^m-1} g_{|J|,|I|}(m-|\BFk|_{I^*}-|\BFd_J)
  \]
  with 
  \[
g_{|J|,|I|}(m-|\BFk|_{I^*}-|\BFd_J) =
\begin{cases}
1+ h_{|J|,|I|}(m-|\BFk|_{I^*}-|\BFd|_J)  & \mbox{ if $m-|\BFk|_{I^*}-|\BFd|_J>2|I|+1$}\\
&\mbox{and $m-|\BFk|_{I^*}-|\BFd|_J$ is odd} \\%and $m>\underline{\nu}_I+1$
\left(\frac{b}{b-1}\right)^{|J|-|I|-1} & \mbox{ if $m - |\BFk|_{I^*}-|\BFd|_J=2|I|+1$}\\
1 & \mbox{ if $m-|\BFk|_{I^*}-|\BFd |_J$ is even,}  
\end{cases}
\]
where
\[
h_{j,i}(\ell) = \frac{b^{j+i-\ell}}{(b-1)^{j-i}} \binom{j-i-1}{\ell-2i}.
\]

(iii) When $m-|\BFk|_{I^*}-|\BFd|_J \le 2|I|$, then $m(\BFk,\BFd,2,J,I)=0$ and therefore we can set $g_{|J|,|I|}(m-|\BFk|_{I^*}-|\BFd|_J) =0$.

\end{proof}

\begin{proof}[of Lemma \ref{lem:Mbs}]
First we observe that
\[
R(b,s) = \left(\frac{b}{2(b-1)}\right)^s P_{\tilde{s},s} ((b-1)/b),
\]
where $P_{m,n}(z)$ is the polynomial defined in Lemma \ref{lem:Ostr}, and recall that $\tilde{s} = \lfloor s/2 \rfloor -1.$
In the notation of \eqref{eq:ostreq}, $z=(b-1)/b$, $z/(z+1) = (b-1)/(2b-1)$, and $1+z = (2b-1)/b$.
Therefore 
\[
R(b,s) = \left(\frac{b}{2(b-1)}\right)^s  \left(\frac{2b-1}{b}\right)^s Pr\left( X  > \frac{b-1}{2b-1}\right),
\]
where $X$ is a Beta rv with parameters $\tilde{s}+1,s-\tilde{s}$.
Now, it is known that a beta distribution  with parameters $a,c$ such that $1<a<c$ has a median no larger than $a/(a+c)$.
%\cite{Ker11}. 
Therefore, if we can show that 
\[
\frac{\tilde{s}+1}{s+1} \le \frac{b-1}{2b-1},
\]
then it means $Pr(X > (b-1)/(2b-1)) \le 1/2$. Now,
\[ 
\frac{b-1}{2b-1} = \frac{1}{2} - \frac{1}{2(2b-1)} \ge \frac{1}{2} - \frac{1}{2(2s-1)}.
\]
On the other hand,
\[ 
\frac{\tilde{s}+1}{s+1} \le \frac{s}{2(s+1)} = \frac{1}{2} - \frac{1}{2(s+1)}.
\] 
Since $2(2s-1) \ge  2(s+1)$ for  $s \ge 2$, we have that
\[
\frac{\tilde{s}+1}{s+1} \le \frac{1}{2} - \frac{1}{2(s+1)} \le \frac{1}{2} -\frac{1}{2(2s-1)} \le \frac{b-1}{2b-1},
\]
as required. So the last step is to show that 
\begin{equation}
\label{eq:yetanotherineq1}
\frac{1}{2} \left(\frac{b}{2(b-1)}\right)^s  \left(\frac{2b-1}{b}\right)^s \le 1.
\end{equation}
The inequality  \eqref{eq:yetanotherineq1} is equivalent to %having
\hchristiane{reword for ref 2; maybe just say in text "which is equivalent to" repeatedly, or omit the first two inequalities; check elsewhere for similar language}
\[
\left(\frac{2b-1}{2b-2}\right)^s \le 2, \mbox{ {to} } %\Leftrightarrow 
s\ln \frac{2b-1}{2b-2} \le \ln 2,  \mbox{ {and finally to} }%\Leftrightarrow 
s \ln (1 + \frac{1}{2(b-1)}) \le \ln 2.
\]
Now, $\ln (1 + \frac{1}{2(b-1)})  \le \ln (1 + \frac{1}{2(s-1)})  \le   \frac{1}{2(s-1)}$, hence it is sufficient to show that
\[
  \frac{1}{2(s-1)} \le \frac{\ln 2}{s} \Leftrightarrow \frac{s}{s-1} \le 2 \ln 2 = 1.3862...
\]
which holds for any $s \ge 4$. For $s=2$, we have that  $R(b,2) = 0.25(b/(b-1))^2$ and since $b \ge s$ we have that $b/(b-1) \le 2$, which implies $0.25(b/(b-1))^2 \le 1$ for all $b \ge 2$. When $s=3$, then  $R(b,3) = 0.125 (b/(b-1))^3 \le 1$.
\end{proof}

\begin{proof}[of Lemma \ref{lem:step1}]
We have that $(m,\mathbf{0}) \in {\cal B}$ is equivalent to assuming $|J|\le m < 2|J|$. We will deal with the case $m = 2|J|-1$ separately, and will first assume $|J| \le m \le 2|J|-3$ and $m$ is odd.

(i) Case where $|J| \le m \le 2|J|-3$ and $m$ is odd: 

In this case, $m^*=m > 2|I|$ if and only if  $|I| < 0.5m$, i.e., $|I| \le  0.5(m-1)$. %where 
%\[
%\tilde{s}_1 = \lceil 0.5(s+j_1) \rceil-1.
%\]
Also, $m^*=m < |I|+|J|$ if and only if  $|I|>m-|J|$. So % if we let $\tilde{g}_1(m,I,\BFk,\BFd) =g_1(m,I,\BFk,\BFd)-1$, then 
$G(m,s,J,\mathbf{0},\mathbf{0})$
%\eqref{eq:BoundCase2}
 is of the form
\begin{align}
G(m,s,J,\mathbf{0},\mathbf{0}) &=\!\!\!\!\!  \sum_{j=0}^{0.5(m-3)}\!\!\! \binom{|J|}{j} 1 + \!\!\!\!\!\!\sum_{j=m-|J|+1}^{0.5(m-3)}  \!\!\binom{|J|}{j} h_j(m)
+ \binom{|J|}{0.5(m-1)} \left(\frac{b}{b-1}\right)^{|I|-0.5(m-1)-1} \notag \\
%\label{eq:BoundCaseB} \\
&=%\left(\frac{1}{2}\right)^s \frac{b^{m}}{b^{m}-1} 
%\left[
\sum_{j=0}^{0.5(m-3)} \binom{|J|}{j}+ \sum_{j = m-|J|+1}^{0.5(m-3)}  \binom{|J|}{j} \binom{|J|-j-1}{m-2j} \frac{(b(b-1))^j}{(b-1)^{|J|}b^{m-|J|}} \notag\\
&+\binom{|J|}{0.5(m-1)}\left(\frac{b}{b-1}\right)^{|J|-0.5(m-1)-1}.
%\right] 
\label{eq:Step1Ineq}
\end{align}

Using Lemma \ref{lem:Sj1incr} with $s=|J|$,  we know that %the term in square brackets in 
\eqref{eq:Step1Ineq} is increasing with $m$, so $G(m,s,J,\mathbf{0},\mathbf{0}) \le G(2|J|-3,s,J,\mathbf{0},\mathbf{0})$ for all $m \le 2|J|-3$
such that $(m,\mathbf{0}) \in {\cal B}$.
%%%%%%%%%%%%%%% GOES AFTER PROOF OF CLAIM
Furthermore, we have that
\begin{align*}
G(2|J|-3,s,J,\mathbf{0},\mathbf{0}) &= \sum_{j=0}^{|J|-3} \binom{|J|}{j} + 
\binom{|J|}{2}\binom{1}{1} \frac{b}{(b-1)} \\
&= 2^{|J|}-1-|J| +\frac{|J|(|J|-1)}{2} \left( \frac{b}{b-1}-1 \right)\\ 
&= 2^{|J|}-1-|J| +\frac{|J|(|J|-1)}{2(b-1)} 
 \le 2^{|J|}-1 - \frac{|J|}{2},
\end{align*}
where the last inequality is obtained by observing that $|J| \le s \le b$.

(ii) Case where $m=2|J|-1$. 

In this case, $m \ge |I|+|J|$ for all $I$ such that $|I| \le |J|-2$.  Also, %recall from Remark \ref{rem:smoins2} that  
for subsets $I$ such that $|I| = |J|-1$, we cannot have $2|I| < m < |I|+|J|$ since in that case
$|I|+|J| - 2|I| = |J|-|I| = 1$.
%
%If $I$ is such that $|I| = |J|-1$ then we cannot have %$2|I| < m^* < |I|+|J|$ since in that case
%$|I|+|J| - 2|I| = |J|-|I| = 1$.
%
%
Therefore, there is no $I$ such that $2|I| < m < |I|+|J|$, and thus % $G(2|J|-1,s,J,\mathbf{0},\mathbf{0})$ is given by
\begin{equation}
\label{eq:G2s-1largest}
G(2|J|-1,s,J,\mathbf{0},\mathbf{0}) = \sum_{j=0}^{|J|-1} \binom{|J|}{j} = 2^{|J|}-1 \ge G(2|J|-3,s,J,\mathbf{0},\mathbf{0}).
\end{equation}

(iii) If $m$ is even with $|J| \le m < 2|J|$, then $G(m,s,J, \mathbf{0},\mathbf{0}) \le \sum_{j=0}^{|J|-1} \binom{|J|}{j} = 2^{|J|}-1$.
\end{proof}

%%%
\begin{proof}[of Lemma \ref{lem:step2}]
We let $j=|J|$ and write
\[
G(m,\BFk) = \sum_{I:|I| \le j-1} g_{j,|I|}(m-|\BFk|_{I^*}).
\]
That is, $G(m,\BFk)=G(m,s,J,\BFk,\mathbf{0})$, i.e., we drop the dependence on $s$, $J$ and $\BFd$. 
%(and also drop the dependence on $\BFd$ in $g(m,J,I,\BFk,\BFd)$ since we have determined that we could set $\BFd=\mathbf{0}$). 
%where
%\[
%g_3(m,\BFk) = \begin{cases}
%1 & \mbox{if } m \ge \bar{\nu}_I \\
%1+ g_1(I,\BFk,\mathbf{0}) & \mbox{ if } \underline{\nu}_I < m \bar{\nu}_I \\
%0 &\mbox{ if }m \le \underline{\nu}_I.
%\end{cases}
%\]

%Claim: $\BFk$ is of the form $(k_1,\ldots,k_s)$ with $k_j$ even. [TODO]

First, we 
%consider the set of subsets $I \subset J$ that contribute a non-zero value to $G(m,\BFk)$. That is, let
%\[
%{\cal G}_{m,\BFk}=\{I: |I| \le j-1, m-|\BFk|_{I^*} \ge 2|I| %+1\}
%\]
define $\iota$ as the size of the largest (strict) subset $I$
of $J$ that contributes a non-zero value to $G(m,\BFk)$. That is, %$j_{min}$ and 
%$\iota$ as the size of the  %smallest and 
%largest %, respectively, 
%subset(s) in ${\cal G}_{m,\BFk}$. That is,
\begin{align*}
\iota &:= \max \{|I| : I \subset J, m-|\BFk|_{I^*} \ge 2|I|+1\}.
%j_{min} &= \min \{j : \exists\,\, I \in {\cal G}_{m,\BFk} \mbox{ with }|I|=j\}.
\end{align*}
Note that $0  \le \iota \le j-1$. \hchristiane{$j=|J|$ was defined at the beginning of the proof, on the previous page. I thought you had suggested to use this instead of $|J|$ and I do find it makes the notation less clunky.} Also, it is useful at this point to mention that our optimal solution $(\tilde{m},\mathbf{0})$ will be such that $\tilde{m} = 2\iota+1$. 

We then define
$
{\cal G}_{m,\BFk} := \{I:0 \le |I| \le \iota\}.
$
%Note that ${\cal G}_{m,\BFk} \subseteq \tilde{\cal G}_{m,\BFk}$ and $\tilde{\cal G}_{m,\BFk}$  may include subsets $I$ such that
%$g_{j,|I|}(m-|\BFk|_{I^*}) = 0$. 
Using this notation we can write
%We also have that 
\begin{equation}
    \label{eq:Gmksum}
G(m,\BFk) = \sum_{I: I \in {\cal G}_{m,\BFk}} g_{j,|I|}(m-|\BFk|_{I^*}).
\end{equation}
%= \sum_{I \in \tilde{\cal G}_{m,\BFk}} g_{j,|I|}(m-|\BFk|_{I^*})$. 
This holds because if $I \notin {\cal G}_{m,\BFk}$, then $m-|\BFk|_{I^*}\le 2|I|$ and thus $g_{j,|I|}(m-|\BFk|_{I^*}) = 0$.

\hchristiane{I saw your suggestion to drop the index $j$, but I'm not convinced yet that this is the best thing to do. Although $j$ does not play a direct role here, Lemma B.1 refers to the definition (11) which does use $j$. So if we drop $j$ we'd have to say something about the mismatch of notation between the statement of Lemma B.1 and (11). I'm not sure having to do this is better than carrying the index through.}

%The following lemma gives us useful properties about the function $g_{j,i}(\ell)$.
%[TODO:refer to lemma imstead]

Next we introduce a definition:\\
%\begin{definition}
%\label{def:dom}
\noindent{\bf Definition:}
For a given $J$ and $I \subset J$,
we say that {\em $(m,\BFk)$ is dominated by $(m',\mathbf{0})$ at $I$} if $g_{j,|I|}(m-|\BFk|_{I^*}) \le g_{j,|I|}(m')$. 
%\end{definition}

Our strategy will be as follows: consider the set ${\cal M} = \{ 1,3,\ldots, 2\iota+1\}$. We claim that 
for each $I \in {{\cal G}}_{m,\BFk}$, there exists $(m',\mathbf{0})$ with $m' \in {\cal M}$ such that 
%such that on each level $j$ and subset $I$ such that $|I| = j$ with $j_{min} \le j \le j_M$, there exists $m_{\ell} \in \%{m_l,m_l+2,\ldots,m_u\}$ such that $g_3(I,\BFk) \le g_3(I,m_{\ell},\mathbf{0})$, i.e., such that 
$(m,\BFk)$ is dominated by $(m',\mathbf{0})$ at $I$. 
\hchristiane{The rest of the paragraph is new.}
In turn, this will allow us to bound each term $g_{j,|I|}(m-|\BFk|_{I^*})$ in \eqref{eq:Gmksum}
by a term of the form $g_{j,|I|}(m')$. We then only need to keep track, for each $m'$, of how many times $g_{j,|I|}(m')$  has been used in this way---something we will do by introducing counting numbers denoted by $\eta(\cdot)$. This strategy is a key intermediate step to get to our end result, which is to show that $G(m,\BFk) \le G(\tilde{m},\mathbf{0})$.

To prove the existence of this $m'$, we define a mapping
${\cal L}_{m,\BFk}: {{\cal G}}_{m,\BFk} \rightarrow {\cal M}$ that will, for a given $m$ and $\BFk$, assign to each subset $I \in {{\cal G}}_{m,\BFk}$, the largest integer $m' \in {\cal M}'$ such that $(m,\BFk)$ is dominated by $(m',\mathbf{0})$ at $I$. The reason why we choose the largest $m' \in {\cal M}$ is that this provides us with the tightest bound on $g_{j,|I|}(m-|\BFk|_{I^*})$, as should be clear from the behavior of the function $g_{j,i}(\ell)$, as described in Lemma \ref{lem:propfji}.

The mapping ${\cal L}_{m,\BFk}(I)$ is defined as follows:
\[
{\cal L}_{m,\BFk}(I) := \begin{cases} 
2\iota+1 & \mbox{ if $m-|\BFk|_{I^*} \ge 2\iota+1$}\\
2\ell+1 & \mbox{  if $2\ell+1\le m-|\BFk|_{I^*} \le 2{\ell}+2$ where $0 \le \ell < \iota$}   \\
1 & \mbox{ if $m-|\BFk|_{I^*}\le 0$.} 
\end{cases}
\]

%where ${\cal L}_{m,\BFk}(I) = 2j'+1$, and $j'$ is the largest integer such that $m-|\BFk|_I \ge 2j'+1$, or, if $m-|\BFk|_I \le 2|I|$, then  ${\cal L}_{m,\BFk}(I) = 2j_{min}+1$ .% the function that maps $I$ to $m'$ in this way. 
%In other words, 
%\[
%{\cal L}_{m,\BFk}(I) = \begin{cases} 
%2j_{max}+1 & \mbox{ if $m-|\BFk|_{I^*} \ge 2j_{max}+1$}\\
%2j'+1 & \mbox{  if $m' \le m-|\BFk|_{I^*} < m'+2$ where %$m'=2j'+1$ and $j' \in \{1,\ldots,j_{max+1}\}$}   \\
%1 & \mbox{ if $m-|\BFk|_{I^*} < 3$.} 
%\end{cases}
%\]

%\begin{claim}
\noindent{\bf Claim:}
For each $I \in \tilde{{\cal G}}_{m,\BFk}$, 
$({\cal L}_{m,\BFk}(I),\mathbf{0})$ dominates $(m,\BFk)$ at $I$.
%\end{claim}

%\begin{proof}
\noindent{\bf Proof:}
We need to show that  $g_{j,|I|}(m-|\BFk|_{I^*}) \le g_{j,|I|}({\cal L}_{m,\BFk}(I))$. We proceed by examining the three possible cases for ${\cal L}_{m,\BFk}(I)$ based on its definition.

(i) Assume $m-|\BFk|_{I^*} \ge 2\iota+1$ and therefore ${\cal L}_{m,\BFk}(I)=2\iota+1$. By definition of $\iota$ we have  
$2\iota+1 > 2i$, and therefore Part 2 of
Lemma \ref{lem:propfji} applies, which implies that
$g_{j,|I|}(2\iota+1) \ge g_{j,|I|}(m-|\BFk|_{I^*})$. % since, as mentioned above, $m-|\BFk|_{I^*} \ge 2\iota+1$.
(ii) We now assume $2\ell+1 < m-|\BFk|_{I^*} \le 2{\ell}+2$ for some $ 0\le \ell < \iota$. In this case, ${\cal L}_{m,\BFk}(I) = 2\ell+1$ and $m-|\BFk|_{I^*}$ is either equal to $2\ell+1$  or to $2\ell+2$. If $m-|\BFk|_{I^*} = 2\ell+1$  then clearly $g_{j,|I|}({\cal L}_{m,\BFk}(I)) \ge g_{j,|I|}(m-|\BFk|_{I^*})$ since in fact these two quantities are equal. If $m-|\BFk|_{I^*}  = 2\ell+2$ then since ${\cal L}_{m,\BFk}(I) =2\ell+1 \ge 1$ is odd we can use Part 1 of Lemma \ref{lem:propfji} to conclude that  $g_{j,|I|}({\cal L}_{m,\BFk}(I) ) \ge g_{j,|I|}(m-|\BFk|_{I^*})$.
(iii) If $m-|\BFk|_{I^*} \le 0$, then $m-|\BFk|_{I^*} \le 2|I|$ since $I \in {{\cal G}}_{m,\BFk}$ implies
$0 \le |I| \le \iota$, and therefore $g_{j,|I|}(m-|\BFk|_{I^*}) =0 \le g_{j,|I|}(1)$. 
\,\,\,$\square$
%\end{proof}

%(Note that $j_{max-i} = j_{max}-i$ and $j_{max-v} = j_{min}$.) 
\hchristiane{If we define them in reverse as you suggest (with ${\cal L}_{m,\BFk}(I) = 2\ell+1$), the main issue that arises is that the interpretation based on probabilities becomes counterintuitive since the range of values for $|\BFk|_{I^*}$ would decrease as $\ell$ increases.
I guess we could define the probabilities in terms of the $m-|\BFk|_{I^*}$ but then it would make the argument with the cumulative probabilities more clumsy, which in turn would make it harder to mention as being obvious (as we now do). At the end, I'm finding there is no ideal way to deal with this, i.e., it seems like however we do it, there is a need at some point to reverse the indices based on $\ell$. For now I'm leaving it in this way because of the probability argument mentioned above, but I'm open to change if I can find a way to do it that doesn't introduce just a different type of clumsiness compared to what we currently have.}

\hchristiane{This paragraph is new}
Now, recall that shortly after stating \eqref{eq:Gmksum}, when we explained our strategy to replace the terms $g_{j,|I|}(m-|\BFk|_{I^*})$ by $g_{j,|I|}(m')$ in \eqref{eq:Gmksum}, we also said we would need counting numbers $\eta(\cdot)$
%in order to get the final result we want, which is to show that $G(m,\BFk) \le G(\tilde{m},\mathbf{0})$, we will need to carefully examine the sum \eqref{eq:Gmksum} once each $g_{j,|I|}(m-|\BFk|_{I^*})$ is replaced by a term of the form $g_{j,|I|}(m')$. In particular, we'll need to know, for fixed $j$ and $|I|$, 
to tell us how many times, for each $m'$, the term $g_{j,|I|}(m')$ has been used in this way. These counting numbers are essential to apply the optimization result involving weighted sums that is given in Lemma \ref{lem:weightedsums}, which is the key to get our final result.  They are defined as follows, for $0 \le i,\ell \le \iota$:
\[
\eta(\ell,i,\BFk) := |\{I \in {{\cal G}}_{m,\BFk}: |I|=i,  {\cal L}_{m,\BFk}(I)=2(\iota-\ell)+1\}|.
\]
\hchristiane{check that it should be   ${\cal G}_{m,\BFk}$ and not $\tilde{{\cal G}}_{m,\BFk}$.}
Note that $\sum_{\ell=0}^{\iota} \eta(\ell,i,\BFk) = \binom{j}{i} $. Also we can think of $p(\ell,i,\BFk) = \eta(\ell,i,\BFk)/\binom{j}{i}$ as the probability that a randomly chosen subset $I$ of $i$ elements from $J$ is such that $|\BFk|_{I^*} \in {\cal R}_{\ell}$, where
\begin{equation}
    \label{eq:Rl}
{\cal R}_{\ell} := \begin{cases}
\{0,1,\ldots,m-(2\iota+1)\} & \mbox{ if } \ell=0\\
\{m-(2\iota+1)+2\ell-1,m-(2\iota+1)+2\ell)\} & \mbox{ if } 1 \le \ell < \iota \\
\{m-2,m-1,\ldots\} & \mbox{ if } \ell=\iota.
\end{cases}
\end{equation}
%\le m_{\ell}$ for $1 \le \ell <\iota$; (ii) $|\BFk|_{I^*} \le m_0$ if $\ell=0$; (iii) $m_{\iota-1} < |\BFk|_{I^*}$ if $\ell=\iota$. 
%Accordingly, we define the cumulative probabilities
%\[
%P(r,i,\BFk) = \sum_{\ell=0}^{r} p(\ell,i,\BFk)
%\]
%for $r=0,\ldots,\iota$.
%where $m'=2j'+1$, $j_{min} \le j ,j'\le j_{max}$. 
%That is, for $r=0,\ldots,\iota-1$, $P(r,i,\BFk) $ is the probability that a randomly chosen subset $I$ of $i$ elements from $J$ is such that $|\BFk|_{I^*}$ is no larger than $m_{r}$, with $P(\iota,i,\BFk)=1$ for each $i$. It should be obvious that 

%\begin{claim}
%\begin{equation}
%\label{claim:decrprob}
%P(r,i,\BFk) \ge P(r,i+1,\BFk) \mbox{ for %$i=0,\ldots,\iota-1$, and $r=0,\ldots,\iota$,}
%\end{equation}
%i.e., the event that the sum of $i$ randomly chosen integers %from the set of non-negative integers in $\BFk_J$ is such %that $|\BFk|_{I^*}  \le m_{r}$ has a probability larger or %equal to the same probability but for a sum of $i+1$ %integers.

%\begin{proof} We have that
%\[
%P(\ell,j+1,\BFk)= \sum_{i=0}^{\ell} p(i,j+1,\BFk)= \sum_{i=0}^{\ell} p(i,j,\BFk)  \pi_{i,\ell,j,\BFk} \le  \sum_{i=0}^{\ell} p(i,j,\BFk) 
%= P(\ell,j,\BFk),
%\]
%where $\pi_{i,\ell,j,\BFk}$ is the probability to  choose an element in $\{k_1,\ldots,k_s\}$ that is $\le m_{\ell}-m_i=2(\ell-i)$ given %$j$ elements in $\{k_1,\ldots,k_s\}$ adding to a value in the range associated with $p(i,j,\BFk)$ have already been chosen.
%\end{proof}

To get the final result, we write:
\begin{align}
G(m,\BFk) & =  \sum_{I: 0 \le |I| \le \iota} g_{j,|I|}(m-|\BFk|_{I^*})  \\
&\le \sum_{i=0}^{\iota}  \sum_{\ell=0}^{\iota} \eta(\ell,i,\BFk) g_{j,i}(2(\iota-\ell)+1)\notag \\
&= \sum_{i=0}^{\iota}  \sum_{\ell=0}^{\iota} p(\ell,i,\BFk) \binom{j}{i} g_{j,i}(2(\iota-\ell)+1)\notag
\end{align}
\begin{align}
&= \sum_{i=0}^{\iota}  \sum_{\ell=0}^{\iota-i} p(\ell,i,\BFk) \binom{j}{i} g_{j,i}(2(\iota-\ell)+1)
\label{eq:cutsum} \\
&\le \sum_{i=0}^{\iota} \binom{j}{i}  g_{j,i}(2\iota+1) = G(\tilde{m},\mathbf{0}),\notag
\end{align}
where $\tilde{m}=2\iota+1$. In the above, the first inequality is obtained by replacing $g_{j,|I|}(m-|\BFk|_{I^*})$ by $g_{j,i}(2(\iota-\ell)+1)$ for each of the $\eta(\ell,i,\BFk)$ pairs $(m,\BFk)$ dominated by $(2(\iota-\ell)+1,\mathbf{0})$ at $I$, where $|I|=i$; the third equality holds because if $\ell>\iota-i$, then 
$
2(\iota-\ell)+1 < 2(\iota-(\iota-i))+1 = 2i+1
$
and therefore $g_{j,i}(2(\iota-\ell)+1)=0$. Similarly, $\ell \le \iota-i$ implies $2(\iota-\ell)+1 \ge 2i+1$ and so 
$g_{j,i}(2(\iota-\ell)+1)>0$ in this case.
The last inequality comes from applying Lemma \ref{lem:weightedsums}, whose conditions hold because:
\begin{enumerate}
\item $\binom{j}{i}g_{j,i}(2(\iota-\ell)+1)$ corresponds to $x_{\ell+1,i+1}$ in Lemma \ref{lem:weightedsums};
\item  Lemma \ref{lem:Sj1incr} together with  \eqref{eq:G2s-1largest} shows the decreasing row-sums condition is satisfied, i.e., $G(2(\iota-\ell)+1,\mathbf{0}) = \sum_i x_{\ell+1,i+1}$ is decreasing with $\ell$; %of Claim \ref{claim:k0} (which we currently are in the process of proving Step 2);
\item The increasing-within-column assumption of  Lemma \ref{lem:weightedsums} is satisfied because the sum \eqref{eq:cutsum} only includes positive values of 
$g_{j,i}(2(\iota-\ell)+1)$ (as shown above), which in turn allows us to invoke Part 2 of Lemma \ref{lem:propfji};
\item $p(\ell,i,\BFk)$ corresponds to $\alpha_{\ell+1,i+1}$ in Lemma \ref{lem:weightedsums};
\item  To see that the $p(\ell,i,\BFk)$'s obey the decreasing-cumulative-sums condition \eqref{eq:decrAij} in Lemma \ref{lem:weightedsums}, we argue that our probabilistic interpretation of the $p(\ell,i,\BFk)$ based on  the sets defined in \eqref{eq:Rl} should make it clear that for $i=0,\ldots,\iota -1$ and $0 \le r \le \iota$,
\[
\sum_{\ell=0}^{r} p(\ell,i,\BFk) \ge \sum_{\ell=0}^{r} p(\ell,i+1,\BFk).
\]
\end{enumerate}
Therefore 
$
G(m,\BFk) \le G(\tilde{m},\mathbf{0}),
$
%where $\tilde{m} = 2j_{max}+1$, 
as required.
\end{proof}
% end of proof of Lemma \ref{lem:step2}

%\section*{Appendix: Technical Lemmas}
%\addcontentsline{toc}{section}{Appendix: Technical Lemmas}
%\label{sec:techlemma}

%%%%%%%%%%%%%%%%%%%%%%%%%%%%%%%%%%%%%%%%%%%
\subsection*{Technical lemmas}

The following result \cite[Lemma2]{Ost04} is
%We start by recalling a result from Ostrovskii \cite{Ost04} 
used to prove intermediate inequalities needed in our analysis.
\hchristiane{Ref 2 wants us to bring more context from that paper since $z \in \mathbb{C}$ and asks what if $z = -1$}
% \begin{lemma}
% \label{lem:Ostr}
% Let $P_{m,n}(z)$ be the polynomial defined by
% \[
% P_{m,n}(z) = \sum_{j=0}^m \binom{n}{j} z^j, \qquad 0 < m < n-1, z \in \mathbb{C}.
% \]
% Then
% \begin{equation}
% \label{eq:ostreq}
% \frac{P_{m,n}(z) }{(1+z)^n \binom{n}{m}(n-m)} = \int_{z/(z+1)}^1 u^{m} (1-u)^{n-m-1}du.
% \end{equation}
% \end{lemma}

% Note that \eqref{eq:ostreq} can be rewritten as
% \[
% P_{m,n}(z) = (1+z)^n \int_{z/(z+1)}^1 \frac{n!}{m!(n-m-1)!} u^{m} (1-u)^{n-m-1}du
% \]
% where the RHS corresponds to the pdf of a beta distribution with parameters $(m+1,n-m).$

{
\begin{lemma}
\label{lem:Ostr}
Let $P_{m,n}(z)$ be the polynomial defined by
\[
P_{m,n}(z) = \sum_{j=0}^m \binom{n}{j} z^j, \qquad 0 < m < n-1, %z \in \mathbb{C}.
\]
Then for $z \neq -1$
%\[
%P_{m,n}(z) = (1+z)^n \int_{z/(z+1)}^1 \frac{n!}{m!(n-m-1)!} u^{m} %(1-u)^{n-m-1}du
%\]
%where the RHS corresponds to the pdf of a beta distribution with %parameters $(m+1,n-m).$
%In particular,
\begin{equation}
\label{eq:ostreq}
\frac{P_{m,n}(z) }{(1+z)^n \binom{n}{m}(n-m)} = \int_{z/(z+1)}^1 u^{m} (1-u)^{n-m-1}du.
\end{equation}
%whenever $z\neq -1$.
\end{lemma}
}

We also need the following identity for {integers $c>a\geq0$, %taken from \cite{Fil07}:
which may be found in \cite[(5.16)]{knuth1989}}
\hchristiane{Ref 1 wants us to cite instead ``Concrete Mathematics'', by Graham, Knuth, Patashnik and Ref 2 says this holds for $a=0$}\hjaspar{I assume you wanted me to track down the reference. Please double check that I entered it in the bibliography correctly.}
\begin{equation}
\label{eq:AlternBinTruncIdent}
\sum_{j=0}^a \binom{c}{j} (-1)^j = \binom{c-1}{a}(-1)^a.
\end{equation}
%%%%%%%%%%%%%%%%%%%%%%%

We now state and prove a number of  technical lemmas that were used within the above proofs.
The first two lemmas are used to prove  Lemma \ref{lem:boundPsi}, and the next three are used for Lemmas \ref{lem:step1} and \ref{lem:step2}.

%%%%%%%%%%%%%%%%%%%%%%%%%%%%%%%%%%%%%%%%%%%%%%%%%
\hchristiane{Ref 2 says this holds for $k=0$: include only if it doesn't mess up the proof; otherwise explain in response why we didn't include it}
\begin{lemma}
\label{lem:Qbks}
For $b \ge s \ge 2$ and ${0 \le }  k < s$,  let
\[
Q(b,k,s) := \sum_{j=0}^{k} (-1)^j \binom{s}{j} (b^{k-j}-1).
\]
Then 
\begin{equation}
\label{eq:TightBoundQbks}
Q(b,k,s) \le \begin{cases}
b^k \left(\frac{b-1}{b}\right)^s & \mbox{ if  $k$ is even}\\
b^k \left(\frac{b-1}{b}\right)^s +\binom{s-1}{k} & \mbox{ if $k>1$ is odd},\\
(b-1) & \mbox{ if $k=1$}.
\end{cases}
\end{equation}
%and for both $k$ odd and even we have
%\begin{equation}
%\label{eq:looseBound}
%%Q(b,k,s) \le (b-1)^k.
%\end{equation}
\end{lemma}

\begin{proof}
{The statement holds trivially for $k=0$. For $k>0$ we apply Lemma \ref{lem:Ostr} and \eqref{eq:AlternBinTruncIdent} to obtain}
\begin{align*}
Q(b,k,s) &= b^k P_{k,s}(-1/b)- \binom{s-1}{k}(-1)^k\\
&= b^k \left(\frac{b-1}{b}\right)^s \int_{-1/(b-1)}^1 u^k (1-u)^{s-k-1} c_{k,s} du - \binom{s-1}{k}(-1)^k \\
&= b^k \left(\frac{b-1}{b}\right)^s \left(\int_{-1/(b-1)}^0 u^k (1-u)^{s-k-1} c_{k,s} du +1\right)  - \binom{s-1}{k}(-1)^k 
\end{align*}
where $c_{k,s}$ is the constant that makes the integrand a beta pdf with parameters $(k+1,s-k)$, i.e., $c_{k,s} = s \binom{s-1}{k}$.

%We need to show that $Q(b,k,s) \le (b-1)^k$ which holds if and only if
%\[
%\frac{Q(b,k,s)}{b^k}=\left(\frac{b-1}{b}\right)^s (\int_{-1/(b-1)}^0 u^k (1-u)^{s-k-1} c_{k,s} du +1) - \frac{\binom{s-1}{k}(-1)^k}{b^k} \le 
%\left(\frac{b-1}{b}\right)^k.
%\]
Now, if $k$ is odd then $\int_{-1/(b-1)}^0 u^k (1-u)^{s-k-1} c_{k,s} du \le 0$ and thus  we get
\[
Q(b,k,s) \le b^k \left(\frac{b-1}{b}\right)^s +\binom{s-1}{k}.
\]
Note also that when $k=1$, $Q(b,k,s) = b-1$ (which is not necessarily bounded from above by $b^k \left(\frac{b-1}{b}\right)^s +\binom{s-1}{k} = b ((b-1)/b)^s+s-1$. It is for this reason we treat the case $k=1$ separately).
If $k$ is even, then
\[
\int_{\frac{-1}{b-1}}^0 \!\!\! u^k (1-u)^{s-k-1} c_{k,s} du \le c_{k,s} \frac{1}{b-1} \frac{1}{(b-1)^k} \left(1+ \frac{1}{b-1}\right)^{s-k-1}
=c_{k,s} \frac{b^{s-k-1}}{(b-1)^s}.
\]
Therefore when $k$ is even
\begin{align*}
Q(b,k,s) &= b^k \left(\frac{b-1}{b}\right)^s  \left( 1+c_{k,s}\frac{b^{s-k-1}}{(b-1)^s} \right)- \binom{s-1}{k} \\
% \le \left(\frac{b-1}{b}\right)^k \\
&=  b^k \left(\frac{b-1}{b}\right)^s +  \left(\frac{s}{b} -1\right) \binom{s-1}{k} \le  b^k \left(\frac{b-1}{b}\right)^s
\end{align*}
since $s\le b$. 
\end{proof}

%[TODO:change the notation]
\begin{lemma}
\label{lem:Mtilde}
Consider a $(0,m,s)$-net in base $b$. %with corresponding quantity $n_b(\BFi)$ %introduced in \cite{WLD20}.
Let $\emptyset \neq J \subset \{1,\ldots,s\}$, 
and $I \subseteq J$, with $I^*=I \cup J^c$.
%let
%\[
%\tilde{n}_{b} (\BFk;\BFd;I^c)= %\sum_{\substack{i_j = d_j, j \in I^c \\ i_{j} %\ge k_{j},j \in I}}
%n_b(\BFi).
%\]
Then the following bounds hold:

\noindent (i) if $m \ge |\BFk|_{I^*} +\BFd_{J}+c|I|+|I^{*c}| $, then
\[
m_b(\BFk,\BFd,c,J,I;P_n)= b^{m-|\BFk|_{I^*} -|\BFd|_{J}-c|I|-|I^{*c}| } (b-1)^{|I^{*c}|};
\]
(ii) if  $|\BFk|_{I^*} + |\BFd|_{J} +c|I|+1 <  m < |\BFk|_{I^*} +\BFd_{J}+c|I|+|I^{*c}|  $
then
\[
m_b(\BFk,\BFd,c,J,I;P_n) \le  b^{m- |\BFk|_{I^*}-|\BFd|_{J}}
\left(\frac{b-1}{b}\right)^{|I^{*c}|}+ \binom{|I^{*c}|-1}{m- |\BFk|_{I^*}-|\BFd|_{J}-c|I|}{\bf 1}_{m-|\BFk|_{I^*}-|\BFd|_{J}-c|I|}
\]
where
\[
{\bf 1}_{x} = \begin{cases}
1 & \mbox{ if } x \mbox{ is odd} \\
0 & \mbox{otherwise}.
\end{cases}
\]

(iii) if $m = |\BFk|_{I^*} + |\BFd|_{J}+c|I| +1 $ then 
\[
m_b(\BFk,\BFd,c,J,I;P_n)  = (b-1).
\]

(iv) if $m \le  |\BFk|_{I^*} + |\BFd|_{J}+c|I| $ then $m_b(\BFk,\BFd,c,J,I;P_n)=0$.
\end{lemma}

\begin{proof}
Using the quantities $n_b(\BFk)$ defined in \cite{WLD20}, their relation to $m_b(\BFk;P_n)$ and the value of the latter for a $(0,m,s)$-net, we write
\begin{align}
m_b(\BFk,\BFd,c,J,I;P_n) &=  \sum_{i_j \ge k_j,j \in J^c; i_j \ge k_j+d_j+c, j \in I} \!\!\!\!\!\! n_b(\BFi_{I^*}: \BFd_{I^{*c}})  \nonumber \\ 
&= \sum_{\BFe \in \{0,1\}^{|I^{*c}|}} (-1)^{|\BFe|} m_b((\BFk_{J^c}:(\BFk+\BFd+2)_I:(\mathbf{d}+\BFe)_{I^{*c}});P_n) \nonumber \\
&=  \sum_{j=0}^{|I^{*c}|} (-1)^j  \binom{|I^{*c}|}{j} \max(b^{m-|\BFk|_{I^*}- |\BFd|_{J}-c|I|- j}-1,0)  \label{eq1}
\end{align}
where $(\BFi_I:\BFd_{I^c})$ represents the vector with $j$th component given by $i_j$ if $j \in I$ and by $d_j$ if $j \notin I$.
If $m-|\BFk|_{I^*}- |\BFd|_{J}-c|I|-|I^{*c}| \ge 0$ then the above sum is given by
\begin{align*}
m_b(\BFk,\BFd,c,J,I;P_n) &=   b^{m- |\BFk|_{I^*} -|\BFd|_{J}-c|I|-|I^{*c}| }\sum_{j=0}^{|I^{*c}|}  (-1)^j  \binom{|I^{*c}|}{j} b^{|I^{*c}|-j}\\
&=   b^{m- |\BFk|_{I^*}-|\BFd|_{j}-c|I|-|I^{*c}| } (b-1)^{|I^{*c}|}.
\end{align*}
If $m-|\BFk|_{I^*} - |\BFd|_{J}-c|I| \le 0$ then the $\max$ inside the sum \eqref{eq1} always yields 0. When $ 1 < m-|\BFk|_{I^*} -
 |\BFd|_{J}-c|I| 
< |I^{*c}|$, then \eqref{eq1} is given by
\begin{align*}
& \sum_{j=0}^{m-|\BFk|_{I^*}- |\BFd|_{J}-c|I|} (-1)^j  \binom{|I^{*c}|}{j} (b^{m-|\BFk|_{I^*}- |\BFd|_{J}-c|I|- j}-1) \\
&\le 
\left(b^{m- |\BFk|_{I^*}-|\BFd|_{J}}\left(\frac{b-1}{b}\right)^{|I^{*c}|} + \binom{|I^{*c}|-1}{m- |\BFk|_{I^*}-|\BFd|_{J}-c|I|}{\bf 1}_{m-
|\BFk|_I-|\BFd|_{I^c}-c|I|}\right),
\end{align*}
\hchristiane{check the $c|I|$ term}
where the last inequality is obtained by applying Lemma \ref{lem:Qbks} with $s=|I^c|$ and $k = m-|\BFk|_{I^*}- |\BFd|_{J}-c|I|$. Finally, when 
$m=|\BFk|_{I^*}+|\BFd|_{J}+c|I|+1$, then \eqref{eq1} is given by
\begin{align*}
& \sum_{j=0}^{1} (-1)^j  \binom{|I^c|}{j} (b^{1- j}-1)  =  b-1. \qed
\end{align*}
 
%[TODO: prove this last inequality, which seems to hold numerically]
\end{proof}

\begin{lemma}
\label{lem:propfji}
The function $g_{j,i}(\ell)$ defined in \eqref{eq:Defgijl} with $0 \le i < j$ and $j \ge 1$ satisfies the following properties.
\begin{enumerate}
%\item  $g_{j,i}(2i+1) \ge g_{j,i}(2i+3)$.
\item For a given $i$, if $\ell \ge 1$ is odd then $g_{j,i}(\ell) \ge g_{j,i}(\ell+1)$.
\item For a given $i$, if $\ell > 2i$ is odd then $g_{j,i}(\ell) \ge g_{j,i}(\ell+r)$ for all $r \ge 0$. 
\end{enumerate}
\end{lemma}

\begin{proof}
%[Proof of Lemma \ref{lem:propfji}]
For (1): if $\ell \ge j+i$ then $g_{j,i}(\ell) = g_{j,i}(\ell+1) = 1$; if $2i<\ell<j+i$, then 
$g_{j,i}(\ell) >1$ while $g_{j,i}(\ell+1)=1$; if $\ell \le 2i$ then $g_{j,i}(\ell)=0$ and since $\ell$ is odd, it means $\ell \le 2i-1$, and thus $\ell+1 \le 2i$, implying that $g_{j,i}(\ell+1)=0$.
For (2): if $\ell \ge j+i$ then $g_{j,i}(\ell+r)=1$ for all $r\ge 0$; if $2i+1<\ell<j+i$, then  the function $g_{j,i}(\ell)$ is increasing as $\ell$ decreases over odd values strictly between $j+i$ and  $2i+1$; this is because when $\ell$ decreases by 2, $h_{j,i}(\ell)$ increases by a factor of at least $2b^2/((j-i-1)(j-i-2))$, which is at least 1 since $j \le s \le b$. Finally, we need to show that $g_{j,i}(2i+1) = (b/(b-1))^{j-i-1} \ge g_{j,i}(2i+3)$, i.e., that   
\[
\left(\frac{b}{b-1}\right)^{j-i-1} \ge 1+ \frac{b^{j+i-(2i+3)}}{(b-1)^{j-i}} \binom{j-i-1}{2i+3-2i} = 1 + \left(\frac{b}{b-1}\right)^{j-i}
 \frac{1}{b^3} \binom{j-i-1}{3}.
\]
Using the bound $((b-1)/b)^j \le (b-j-1)/(b-1)$ shown in the proof of Lemma \ref{lem:Sj1incr}, we have that the above holds if
\begin{align*}
\frac{b-1}{b} &\ge  \frac{(j-i-1)(j-i-2)(j-i-3)}{6b^3}+ \frac{b-(j-i)-1}{b-1} \\
\Leftrightarrow \frac{(b-1)^2-b(b-j+i-1)}{b(b-1)} &\ge \frac{(j-i-1)(j-i-2)(j-i-3)}{6b^3} \\
\Leftrightarrow 6b^2(1+b(j-i-1)) & \ge (b-1)(j-i-1)(j-i-2)(j-i-3), 
\end{align*}
which is clearly true since $j\le s \le b$ and $i \ge 0$ and therefore $6b^3(j-i-1) \ge (b-1)(j-i-1)(j-i-2)(j-i-3)$.
\end{proof}

%%%%%%%%%%%%%%%%%%%%%%%%%%%%%%%%%%%%%%%%%%%%%%%%%%%%%%
\begin{lemma}
\label{lem:Sj1incr}
%[This lemma shows that $G(m,s,\mathbf{0},\mathbf{0})$ is increasing as $m$ increases, for $0<m< 2s$.]
Let $s \ge 3$ and $b \ge s$.
Let $m$ be odd with $1 \le m \le 2s-3 $, and consider the function
\[
G(m,s) = %\frac{1}{2^s}\frac{b^{s+j_1}}{b^{s+j_1}-1} \left(
\sum_{j=0}^{0.5(m-3)} \binom{s}{j} + \sum_{j=\max(0,m-s+1)}^{0.5(m-3)} \binom{s}{j} h_{s,j}(m) + \binom{s}{0.5(m-1)}\left(\frac{b}{b-1}\right)^{s-0.5(m-1)-1},
%\right).
\]
where $h_{s,j}(m)$ is as defined in \eqref{eq:defhijl}, i.e.,
\[
h_{s,j}(m) = \binom{s-j-1}{m-2j}
\frac{b^{s+j-m}}{(b-1)^{s-m}}.
\]
Then $G(m,s) \ge G(m-2,s)$ for $m \ge 3$ odd.
That is, $G(m,s)$ is decreasing over the odd integers from $2s-3$ down to 3.
\end{lemma}

\begin{proof}
First, we compute 
\begin{align*}
&\left(\frac{b}{b-1}\right)^{s-0.5(m-1)-1} - h_{s,0.5(m-1)}(m)\\
 &= \left(\frac{b}{b-1}\right)^{s-0.5(m-1)-1} -
\binom{s-0.5(m-1)-1}{1} \frac{b^{0.5(m-1)+s-m}}{(b-1)^{s-0.5(m-1)}} \\
&= \left(\frac{b}{b-1}\right)^{s-0.5(m-1)} \frac{b-1-(s-0.5(m-1)-1)}{b}\\
&=\left(\frac{b}{b-1}\right)^{s-0.5(m-1)} \frac{b-(s-0.5(m-1))}{b}.
\end{align*}
Using this, we can write
\begin{align}
G(m,s) &= %\frac{1}{2^s}\frac{b^{s+j_1}}{b^{s+j_1}-1} \left(
\sum_{j=0}^{0.5(m-3)} \binom{s}{j} + \sum_{j=\max(0,m-s+1)}^{0.5(m-1)} \binom{s}{j} h_{s,j}(m) \nonumber\\
&+ \binom{s}{0.5(m-1)}\left(\frac{b}{b-1}\right)^{s-0.5(m-1)} \frac{0.5(m-1)+b-s}{b}. \label{eq:Gmsready}
%\right).
\end{align}
Next, we show that for $2 \le j \le 0.5(m-1)$: 
\begin{equation}
\label{eq:hjdecr}
\binom{s}{j}h_{s,j}(m) \ge \binom{s}{j-2}h_{s,j-2}(m-2).
\end{equation}
\begin{align*}
&\binom{s}{j}h_{s,j}(m) -\binom{s}{j-2}h_{s,j-2}(m-2)\\
 &= \binom{s}{j} \binom{s-j-1}{m-2j} \frac{b^{s+j-m}}{(b-1)^{s-j}} -
\binom{s}{j-2}\binom{s-(j-2)-1}{m-2-2(j-2)} \frac{b^{s+j-2-(m-2)}}{(b-1)^{s-(j-2)}} \\
&=\binom{s}{j} \binom{s-j-1}{m-2j} \frac{b^{s+j-m}}{(b-1)^{s-j}}-\binom{s}{j-2}\binom{s-j+1}{m-2j+2}\frac{b^{s+j-m}}{(b-1)^{s-j+2}} \\
&=\binom{s}{j}\binom{s-j-1}{m-2j}\frac{b^{s+j-m}}{(b-1)^{s-j}} \left(1-\frac{j(j-1)}{(s-j+2)(s-j+1)} \frac{(s-j+1)(s-j)}{(m-2j+2)(m-2j+1)} \frac{1}{(b-1)^2}\right).
\end{align*}
Hence to prove \eqref{eq:hjdecr}, we need to show that
\[
1 \ge \frac{j(j-1)}{(s-j+2)(s-j+1)} \frac{(s-j+1)(s-j)}{(m-2j+2)(m-2j+1)} \frac{1}{(b-1)^2},
\]
which holds because
\begin{align*}
&\frac{j(j-1)}{(s-j+2)(s-j+1)} \frac{(s-j+1)(s-j)}{(m-2j+2)(m-2j+1)} \frac{1}{(b-1)^2} \le \frac{j(j-1)}{6} \frac{1}{(b-1)^2}\\
&\le \frac{(s-2)(s-3)}{6(b-1)^2} \le \frac{1}{6} \le 1,
\end{align*}
since $j \le 0.5(m-1) \le s-2$ and $b \ge s$.

Using \eqref{eq:Gmsready}, $G(m,s) \ge G(m-2,s)$ can be shown to hold if
\begin{align}
&\sum_{j=0}^{0.5(m-3)} \binom{s}{j} + \sum_{j=\max(0,m-s+1)}^{0.5(m-1)} \binom{s}{j} h_{s,j}(m) \notag\\ 
&+ \binom{s}{0.5(m-1)}\left(\frac{b}{b-1}\right)^{s-0.5(m-1)} \frac{0.5(m-1)+b-s}{b} \notag\\
\ge &
\sum_{j=0}^{0.5(m-5)} \binom{s}{j} + \sum_{j=\max(0,m-2-s+1)}^{0.5(m-3)} \binom{s}{j} h_{s,j}(m-2)\notag\\ &+ \binom{s}{0.5(m-3)}\left(\frac{b}{b-1}\right)^{s-0.5(m-3)} \frac{0.5(m-3)+b-s}{b}.\label{eq:above25}
\end{align}
In turn, using \eqref{eq:hjdecr}, we know that:
\[
\sum_{j=\max(0,m-s+1)}^{0.5(m-1)} \binom{s}{j} h_{s,j}(m) 
\ge \sum_{j=\max(0,m-2-s+1)}^{0.5(m-5)} \binom{s}{j} h_{s,j}(m-2)
\]
and therefore to show \eqref{eq:above25} it is sufficient to show that
\begin{align}
&\binom{s}{0.5(m-3)} +  \binom{s}{0.5(m-1)}\left(\frac{b}{b-1}\right)^{s-0.5(m-1)} \frac{0.5(m-1)+b-s}{b} \notag \\
\ge & \binom{s}{0.5(m-3)}h_{s,0.5(m-3)}(m-2) \notag \\
+&  \binom{s}{0.5(m-3)}\left(\frac{b}{b-1}\right)^{s-0.5(m-3)} \frac{0.5(m-3)+b-s}{b},  \label{eqGmscas1}
\end{align}
where
\begin{align*}
h_{s,0.5(m-3)}(m-2) &= \binom{s-0.5(m-3)-1}{m-2-(m-3)}\frac{b^{s+0.5(m-3)-(m-2)}}{(b-1)^{s-0.5(m-3)}}\\
&=(s-0.5(m-1) )\left( \frac{b}{b-1} \right)^{s-0.5(m-1)} \frac{1}{b-1}.
\end{align*}

The following inequality will be helpful in this proof:

\begin{claim}
For $b  \ge 2$ and $j \ge 1$, we have that 
\begin{equation}
\label{claim:petitesineq}
\left(\frac{b}{b-1}\right)^j \le \frac{b-1}{b-(j+1)}.
\end{equation}
\end{claim}

\begin{proof}
The inequality is equivalent to having
\[
(b-1)^{j+1} \ge b^j (b-(j+1)).
\]
{Applying the mean value theorem to $f(x) = x^{j+1}$ and noticing $f'(x)$ is monotone increasing for $x \ge 0$ , we get that
$f(b)-f(b-1) = f'(\xi) \le f'(b)$ for some $\xi \in (b-1,b)$ and thus $(b-1)^{j+1} \ge b^{j+1} - (j+1)b^j$.}
%From the 
%binomial theorem, we have that $(b-1)^{j+1} \ge b^{j+1} - (j+1)b^j$, as required.
\hchristiane{Ref 1 says bin thm doesn't work because alternating terms are not decreasing in absolute value: says we should use mean value theorem instead for the fct $f(x) = x^{j+1}$ and the monotonicity of its derivative to get $f(b)-f(b-1) = f'(\xi) \le f'(b)$ for some $\xi \in (b-1,b)$}
\end{proof}

Going back to  our goal of proving \eqref{eqGmscas1}, it is sufficient to show that
\begin{align}
&1+\left(\frac{b}{b-1}\right)^{s-0.5(m-1)}  \frac{s-0.5(m-3)}{0.5(m-1)} \frac{0.5(m-1)+b-s}{b}\notag\\
& \ge \left(\frac{b}{b-1}\right)^{s-0.5(m-1)} \left(\frac{s-0.5(m-1)}{b-1} + \frac{b}{b-1}\frac{0.5(m-3)+b-s)}{b} \right) \notag\\
\Leftrightarrow& \left(\frac{b-1}{b}\right)^{s-0.5(m-1)} +\frac{s-0.5(m-3)}{0.5(m-1)} \frac{0.5(m-1)+b-s}{b}\notag \\
&\ge  \left( \frac{s-0.5(m-1)}{b-1} + \frac{0.5(m-3)+b-s}{b-1} \right). \label{theo}
\end{align}
Using \eqref{claim:petitesineq} to simplify the LHS of \eqref{theo}, we see that \eqref{theo} holds if
\begin{align*}
&\frac{b-(s-0.5(m-1))-1}{b-1}+\frac{s-0.5(m-3)}{0.5(m-1)} \frac{0.5(m-1)+b-s}{b} \\
&\ge 
\left( \frac{s-0.5(m-1)}{b-1} + \frac{0.5(m-3)+b-s}{b-1} \right) \\
\Leftrightarrow&
\frac{s-0.5(m-3)}{0.5(m-1)} \frac{0.5(m-1)+b-s}{b} \ge  \frac{s-0.5(m-1)}{b-1}\\
\Leftrightarrow &
\frac{b-1}{b} \ge \frac{s-0.5(m-1)}{s-0.5(m-3)} \frac{0.5(m-1)}{0.5(m-1)+b-s},
\end{align*}
which holds because $s-0.5(m-3) \le s \le b$ and $\frac{b-1}{b} \ge \frac{b-1-x}{b-x}$ if $x \ge 0$.
\end{proof}

\begin{definition}
We denote by ${\cal A}_w$ the set of $w \times \ell$ \emph{weight matrices} $A$ with entries $\alpha_{i,j} \ge 0$ that satisfy the following two  conditions: %.
%Next consider  a  matrix $A$ of weights .
%(Note that the weights $\alpha_{i,j}$ could be positive for indices $(i,j)$ such %that $x_{i,j}=0$.)  
%Assume $A$ satisfies the two following conditions: first, %
1) $\sum_{i=1}^w \alpha_{i,j} = 1$
for all $j=1,\ldots,\ell$; %Let $A$ be the matrix whose element on the $j$th row and $i$th column is $\alpha_{j,i}$. 
%That is, the weights $\alpha_{j,i}$ for $j=1,\ldots,w$ take a convex combination of the elements on the $i$th column of $X$, for $i=1,\ldots,\ell$. 
2) the weights $\alpha_{i,j}$ obey a {\em decreasing-cumulative-sums} condition as follows: for $1 \le i\le w,1 \le j \le \ell$, let 
\begin{equation}
\label{eq:decrAij}
A_{i,j} = \sum_{k=1}^i \alpha_{k,j}.
\end{equation}
Then $A_{i,j} \ge A_{i,j+1}$ for each $i=1,\ldots,w$ and $j=1,\ldots,\ell-1$ (when $i=w$ we have $A_{w,j}=1$ for all $j$).
This means the weights on the first row are decreasing from left to right; the partial sums of the two first rows are decreasing from left to right, etc.  
\end{definition}

%%%%%%%%%%%%%%%%%%%%%%%%%%%%%%%%%%%%%%%%%
\begin{lemma}
\label{lem:weightedsums}
Let $X$ be a $w \times \ell$ matrix  with $\ell \ge w$ and entries $x_{i,j} \ge 0$ and of the form
\[
X = \begin{bmatrix}
x_{1,1} &\cdots  & x_{1,\ell-w+1} & \cdots && x_{1,\ell-1} & x_{1,\ell} \\
x_{2,1} &\cdots  & x_{2,\ell-w+1}  &  \cdots &  & x_{2,\ell-1} & 0 \\
\vdots & &  \vdots& & \iddots  & \iddots& \\
x_{w-1,1} & \cdots & x_{w-1,\ell-w+1} & x_{w-1,\ell-w+2} & 0 & \cdots & \\
x_{w,1} & \cdots & x_{w,\ell-w+1} & 0 & \cdots & & 0 \\
\end{bmatrix}.
\]
that is $x_{i,j}>0$ if and only if $i +j \le \ell+1$, for $1 \le i \le w, 1 \le j \le \ell$. 
We assume $X$ satisfies the following two conditions: first, 
\[
 \sum_{j=1}^{\ell} x_{1,j} \ge \sum_{j=1}^{\ell-1} x_{2,j} \ge \ldots \ge \sum_{j=1}^{\ell-w+1} x_{w,j}
\]
(we refer to this as the {\em decreasing-row-sums} condition) and second,
\begin{equation}
\label{eq:incrpercolumn}
x_{1,j} \le x_{2,j} \le \ldots \le x_{\min(w,\ell-i+1),j}, \qquad j=1,\ldots,\ell 
\end{equation}
(we refer to this as the {\em increasing-within-column} condition).
%and $x_{j,i} = 0$ when $i > \ell-j+1$. Visually, the $x_{j,i}$'s define an $\ell \times w$ matrix
%\[
%X = \begin{bmatrix}
%x_{1,1} &\cdots  & x_{1,\ell-w+1} & \cdots && x_{1,\ell-1} & x_{1,\ell} \\
%x_{2,1} &\cdots  & x_{2,\ell-w+1}  &  \cdots &  & x_{2,\ell-1} & 0 \\
%\vdots & &  \vdots& & \iddots  & \iddots& \\
%x_{w-1,1} & \cdots & x_{w-1,\ell-w+1} & x_{w-1,\ell-w+2} & 0 & \cdots & \\
%x_{w,1} & \cdots & x_{w,\ell-w+1} & 0 & \cdots & & 0 \\
%\end{bmatrix}.
%\]
%The conditions on the $x_{i,j}$ are that the row sums are decreasing but on each column $i$, the elements  from row 1 to row $\ell-i+1$ are increasing.

Let $A$ be a weight matrix in ${\cal A}_w$ and let
\[
\| A \circ X \|_1 = \sum_{j=1}^{\ell} \alpha_{1,j} x_{1,j} + \sum_{j=1}^{\ell-1} \alpha_{2,j} x_{2,j} + \ldots + \sum_{j=1}^{\ell-w+1} \alpha_{w,j} x_{w,j}.
\]

Then %under these conditions on $X$ and $A$, 
for any $A \in {\cal A}_w$ 
\begin{equation}
\label{eq:boundWeightedSums}
\| A \circ X \|_1  \le \sum_{j=1}^{\ell} x_{1,j}.
\end{equation}
That is, the weight matrix  $A \in {\cal A}_w$ that maximize the LHS of \eqref{eq:boundWeightedSums} is the one with 1's on the first row and 0's elsewhere.
\end{lemma}

\begin{proof}

First, note that $A \in {\cal A}_w$ implies that cumulative sums from the last row up are increasing, i.e.,
for $R_{i,j} = \sum_{k=i}^w \alpha_{k,j}$, we have $R_{i,j} \le R_{i,j+1}$ for $j=1,\ldots,\ell-1$.

We proceed by induction on $w \ge 2$.

If $w=2$, then it suffices to show that  for $A \in {\cal A}_2$, we have that
\[
\sum_{j=1}^{\ell}\alpha_{1,j} x_{1,j} +   \sum_{j=1}^{\ell-1} (1-\alpha_{1,j})x_{2,j} \le \sum_{i=1}^{\ell} x_{1,j}
\Leftrightarrow
\sum_{j=1}^{\ell-1} (1-\alpha_{1,j})x_{2,j} \le \sum_{j=1}^{\ell} (1-\alpha_{1,j}) x_{1,j},
\]
or, equivalently, that
\[
\sum_{j=1}^{\ell-1} (1-\alpha_{1,j}) (x_{2,j}-x_{1,j}) \le (1-\alpha_{1,\ell}) x_{1,\ell}.
\]
Now, we know that
\[
\sum_{j=1}^{\ell} x_{1,j} \ge \sum_{j=1}^{\ell-1} x_{2,j} \Leftrightarrow \sum_{j=1}^{\ell-1} (x_{2,j}-x_{1,j}) \le x_{1,\ell}
\]
with $x_{2,j} - x_{1,j} \ge 0$. Therefore
\[
\sum_{j=1}^{\ell-1} (1-\alpha_{1,j})(x_{2,j}-x_{1,j}) \le (1-\alpha_{1,\ell}) \sum_{j=1}^{\ell-1}  (x_{2,j}-x_{1,j})
\le (1-\alpha_{1,\ell})  x_{1,\ell},
\]
where the first inequality holds because the $\alpha_{1,j}$'s are decreasing.

Now assume the statement holds for $w-1\ge 2$. 
%Consider the case where we have $w$ sums with 
%\[
%\sum_{j=1}^{\ell} x_{1,j} \ge \sum_{j=1}^{\ell-1} x_{2,j} \ge \ldots \ge \sum_{j=1}^{\ell-w+1} x_{w,j}
%\%]
%and weights $A\in {\cal A}_w$. 
First we create a new weight matrix  $\tilde{A}$ by merging the two last rows into the second-to-last one and setting the last one to zero, i.e., we define  
$\tilde{\alpha}_{w-1,j}$  as
\begin{align*}
\tilde{\alpha}_{w-1,j} &= \alpha_{w-1,j}+\alpha_{w,j} \qquad j=1,\ldots,\ell \\
\tilde{\alpha}_{w,j} &=0 \qquad j=1,\ldots,\ell \\
\tilde{\alpha}_{i,j} &= \alpha_{i,j}, i=1,\ldots,w-2,j=1,\ldots,\ell.
\end{align*}
%In this way, $\tilde{x}_{w-1,\ell-w+2} = x_{w-1,\ell-w+2} \ge
%\sum_{i=1}^{\ell-w+1} 
With this change, we claim that $\tilde{A} \in {\cal A}_w$. Indeed:
\begin{enumerate}
\item $\tilde{\alpha}_{i,j} \ge 0$
\item $\sum_{i=1}^w \tilde{\alpha}_{i,j} = \sum_{i=1}^{w-2} \alpha_{i,j} + (\alpha_{w-1,j}+\alpha_{w,j}) + 0 = 1$.
\item $\tilde{A}_{i,j} = A_{i,j}$ for $i=1,\ldots,w-2$ and $\tilde{A}_{w-1,j} = A_{w,j} = 1$ for $j=1,\ldots,\ell$.
\end{enumerate}

Next, we show that
\begin{equation}
\label{eq:FirstIneqInductionw}
\| \tilde{A} \circ X \|_1  \ge \| A \circ X\|_1.
\end{equation}
%where
%\begin{align*}
%\tilde{A} \circ X \|_1 &= \sum_{j=1}^w \sum_{i=1}^{\ell} \tilde{\alpha}_{j,i} x_{j,i} \\
%A \circ X \|_1 &= \sum_{j=1}^w \sum_{i=1}^{\ell} \alpha_{j,i} x_{j,i}.
%\end{align*}
Since $\alpha_{i,j} = \tilde{\alpha}_{i,j}$ for $i< w-1$, then \eqref{eq:FirstIneqInductionw} holds if and only if
\begin{align*}
&\sum_{j=1}^{\ell-w+2} (\alpha_{w-1,j}+\alpha_{w,j}) x_{w-1,j} \ge \sum_{j=1}^{\ell-w+2} \alpha_{w-1,j} x_{w-1,j}
+ \sum_{j=1}^{\ell-w+1} \alpha_{w,j} x_{w,j} \\
\Leftrightarrow &\sum_{j=1}^{\ell-w+1} \alpha_{w,j} x_{w-1,j} + \alpha_{w,\ell-w+2} x_{w-1,\ell-w+2} \ge 
\sum_{j=1}^{\ell-w+1} \alpha_{w,j} x_{w,j}  \\
\Leftrightarrow & \sum_{j=1}^{\ell-w+1} \alpha_{w,j} (x_{w,j}-x_{w-1,j}) \le \alpha_{w,\ell-w+2} x_{w-1,\ell-w+2} .
\end{align*}
By the decreasing-row-sum\hjaspar{changed from ``our"} assumption on the $x_{i,j}$'s we know that
\[
0 \le \sum_{j=1}^{\ell-w+1} (x_{w,j}-x_{w-1,j}) \le x_{w-1,\ell-w+2}
\]
and by assumption that $A\in {\cal A}_w$ we have that $\alpha_{w,1} \le \alpha_{w,2}\le \ldots \le \alpha_{w,\ell}$.
Therefore
\begin{align*}
&\sum_{j=1}^{\ell-w+1} \alpha_{w,j} (x_{w,j}-x_{w-1,j}) \le \alpha_{w,\ell-w+1} \sum_{j=1}^{\ell-w+1} (x_{w,j}-x_{w-1,j}) \\
&\le \alpha_{w,\ell-w+1}   x_{w-1,\ell-w+2} \le \alpha_{w,\ell-w+2}   x_{w-1,\ell-w+2}, 
\end{align*}
as required to show that \eqref{eq:FirstIneqInductionw} holds.

Next, to use the induction hypothesis, we observe that  $\tilde{\alpha}_{w,j} = 0$ implies we can essentially ignore the $x_{w,j}$'s. More formally, let $\tilde{A}_{w-1}$ be the matrix formed by the first $w-1$ rows of $\tilde{A}$ and similarly for $X_{w-1}$. %Moreover, let $\tilde{\Balpha}_{w-1}$ be the set of weights $\{w_{j,i},j=1,\ldots,w-1,i=1,\ldots,\ell\}$. 
Then 
$\tilde{A}_{w-1} \in {\cal A}_{w-1}$, since
\begin{enumerate}
\item $\tilde{\alpha}_{i,j} \ge 0$ for $i=1,\ldots,w-1$, $j=1,\ldots,\ell$
\item $\sum_{i=1}^{w-1} \tilde{\alpha}_{i,j} = \sum_{i=1}^w \alpha_{i,j} = 1$ for $j=1,\ldots, \ell$
\item $\tilde{A}_{i,j} \ge \tilde{A}_{i,j+1}$ as verified earlier (and note that $\tilde{\alpha}_{w-1,1} \le \ldots \le \tilde{\alpha}_{w-1,\ell}$
by assumption that $A \in {\cal A}_w$ and since $\tilde{\alpha}_{w-1,j} = R_{w-1,j}$.)
\end{enumerate}
By applying the induction hypothesis, we obtain
\[
\| \tilde{A}_{w-1} \circ X_{w-1} \|_1 \le \sum_{j=1}^{\ell} x_{1,j}
\]
and since
$
\| A \circ X \|_1 \le \| \tilde{A} \circ X \|_1 = \| \tilde{A}_{w-1} \circ X_{w-1} \|,
$
this proves the result.
%still have weights that add up to 1 (i.e., \sum_{j=1}^{w-2} \alpha_{j,i} + \tilde{\alpha}_{w-1,i} = 1$),
%increasing values ($\tilde{x}_{w-1,i} \ge x_{w-2,i}$ for $i=1,\ldots,\ell-w+2$), and 
%decreasing sums ($\sum_{i=1}^{\ell-(w-2)+1} \alpha_{w-2,i} x_{w-2,i} \ge \sum_{i=1}^{\ell-(w-1)+1} \alpha_{w-1,i} x_{w-1,i} 
%= \sum_{i=1}^{\ell-(w-1)+1} \tilde{\alpha}_{w-1,i} \tilde{x}_{w-1,i}.$ 
\end{proof}

\end{document}